\documentclass{article}
\usepackage{amsfonts}
\usepackage{amsthm}
\usepackage{latexsym}
\usepackage{amssymb}
\usepackage{amstext}
\usepackage{amsmath}
\usepackage{supertabular}
\usepackage{epsfig}

\usepackage{graphicx}
\usepackage{amssymb}
\usepackage{epstopdf}
\theoremstyle{plain}
\newtheorem{thm}{Theorem}
\newtheorem{prop}{Proposition}
\newtheorem{cor}{Corollary}
\newtheorem{lem}{Lemma}

\theoremstyle{definition}

\newtheorem{defn}{Definition}
\newtheorem{ex}{Example}

\theoremstyle{remark}

\def\G{\overline{G}}
\def\Kpi{\overline{K}_{2,n}(\pi)}
\def\Ksig{\overline{K}_{2,n}(\sigma)}
\def\K2n{K_{2,n}}

\def\1{\bigskip\noindent}
\def\2{\bigskip}

\begin{document}

\begin{center} {\bf \Large The Homomorphism Poset of  $K_{2,n}$}  \end{center}
\bigskip

\begin{center}
\textsc{Sally Cockburn}

\textsc{Department of Mathematics}

\textsc{Hamilton College, Clinton, NY 13323}

\textsc {\sl scockbur@hamilton.edu}
\end{center}

\begin{center}
\textsc{Yonghyun Song}

\textsc{Department of Mathematics}

\textsc{Hamilton College, Clinton, NY 13323}

\textsc {\sl ysong@hamilton.edu}
\end{center}

\begin{center} {\bf Draft as of \today} \end{center}

\begin{abstract}
A \emph{geometric graph} $\overline{G}$ is a simple graph $G$ together with a straight line drawing of  $G$ in the plane with the vertices in general position.  Two geometric realizations of a simple graph are \emph{geo-isomorphic} if there is a vertex bijection between them that preserves vertex adjacencies and non-adjacencies, as well as edge crossings and non-crossings. 
A natural extension of graph homomorphisms, \emph{geo-homomorphisms},  can be used to define a partial order on the set of geo-isomorphism classes.  
In this paper, the \emph{homomorphism poset} of $K_{2,n}$ is determined by establishing a correspondence between realizations of $K_{2,n}$ and permutations of $S_n$, in which edge crossings correspond to inversions.  
Through this correspondence, geo-isomorphism defines an equivalence relation on $S_n$, which we call {\it geo-equivalence}.
The number of geo-equivalence classes is provided for all $n \leq 9$.  The modular decomposition tree of permutation graphs is used to prove some results on the size of geo-equivalence classes.  
A complete list of geo-equivalence classes and a Hasse diagram of the poset structure are given for  $n  \leq 5$.  
\end{abstract}

\section{Introduction}

A geometric graph $\overline{G}$ is a simple graph $G = \big( V(G), E(G) \big)$ together with a straight line drawing of  $G$  in the plane with vertices in general position, so that no three vertices are collinear and no three edges cross at a single point. (Such a drawing is also called a rectilinear drawing of $G$.)   Any simple graph will have uncountably many geometric realizations, but we identify those that have the same pattern of edge crossings.  This is formalized by extending the definition of graph isomorphism in a natural way to geometric graphs.

\begin{defn}
Let $\overline{G}, \overline{H}$ be geometric realizations of simple graphs $G, H$ respectively. 
A {\it geo-isomorphism} $f:\overline{G} \to \overline{H}$ is a vertex bijection $f:V(G) \to V(H)$ such that for all $u, v, x, y \in V(G)$,
	\begin{enumerate}
	\item $uv \in E(G)$ if and only if $f(u)f(v) \in E(H)$, and
	\item $xy$ crosses $uv$ in $\overline{G}$ if and only if $f(x)f(y)$ crosses $f(u)f(v)$ in $\overline{H}$.
	\end{enumerate}	
\end{defn}

If there exists a geo-isomorphism $f:\overline{G} \to \overline{H}$, we write $\G \cong \overline{H}$.  Geo-isomorphism clearly defines an equivalence relation on the set of all geometric realizations of a simple graph $G$.  A natural impulse is to classify all geometric realizations of a given graph into geo-isomorphism classes. In \cite{BCDM}, Boutin, Cockburn, Dean and Margea have done this for paths $P_n$, cycles $C_n$ and cliques $K_n$, for $n \leq 6$.

\medskip

Graph homomorphisms are a relaxation of graph isomorphisms; they preserve adjacency, but not non-adjacency.  First introduced almost half a century ago, they are the subject of growing interest in graph theory circles. For an excellent survey of this subject, see \cite{HN}.  In \cite{BC}, Boutin and Cockburn extended the definition of graph homomorphisms to geometric graphs.
	
\begin{defn}
Let $\overline{G}, \overline{H}$ be geometric realizations of simple graphs $G, H$ respectively. 	
 A {\it geo-homomorphism} $f:\overline{G} \to \overline{H}$ is a vertex function $f:V(G) \to V(H)$ such that for all $u, v, x, y \in V(G)$,
	\begin{enumerate}
	\item if $uv \in E(G)$, then  $f(u)f(v) \in E(H)$, and
	\item if $xy$ crosses $uv$ in $\overline{G}$, then $f(x)f(y)$ crosses $f(u)f(v)$ in $\overline{H}$.
	\end{enumerate}
\end{defn}

Concentrating on vertex functions that satisfy (1) of Definition 1 and (2) of Definition 2 allows us to define a relation on the set of geometric realizations of a given graph.  

\begin{defn}
Let $\overline{G}$ and $ \widehat{G}$ be geometric realizations of a simple graph $G$. Then set 
$
\overline{G} \preceq \widehat{G}
$
if and only if  there exists a geo-homomorphism 
$f:\overline{G} \to \widehat{G}$
whose underlying map  $f: G \to G $  is a graph isomorphism.
\end{defn}

It is not difficult to see that this relation is both reflexive and transitive. 
To show that it is anti-symmetric, observe that if 
$f:\overline{G} \to \widehat{G}$ is a geo-homomorphism that is also a graph isomorphism, 
then  the total number of edge crossings in $\widehat{G}$ must be at least as big as the total number of edge crossings in $\G$.  
Hence, if we also have $\widehat{G} \preceq \overline{G}$, 
then $\G$ and $\widehat{G}$ must have the same total number of edge crossings. 
This implies that$f$ is in fact a geo-isomorphism, so the relation defined above is in fact  a partial order.

\begin{defn}
The {\it homomorphism poset} $\mathcal{G}$ of a simple graph $G$ is the set of geo-isomorphism classes of its realizations partially ordered by the relation above.
\end{defn}

\begin{figure}[htbp] 
   \centering
   \includegraphics[width=2.5in]{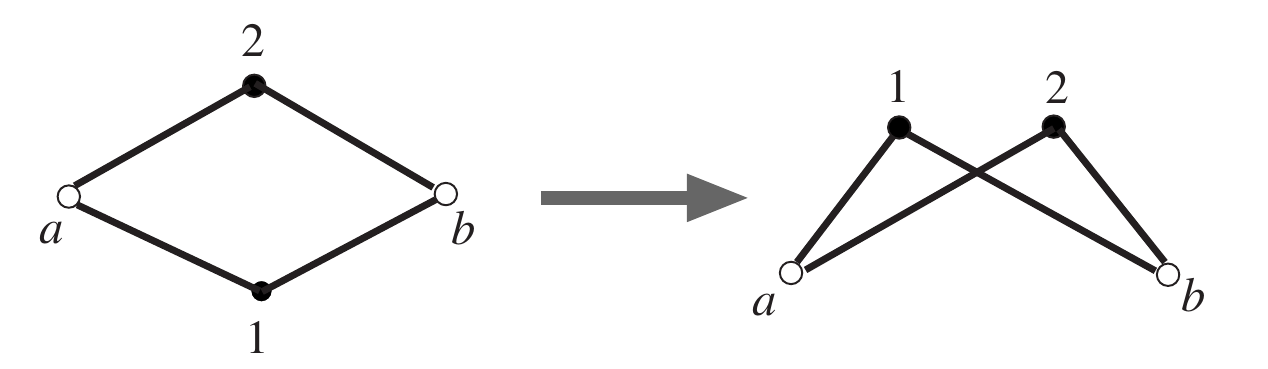} 
   \caption{The homomorphism poset $\mathcal{K}_{2,2}$.}
   \label{fig:K22}
\end{figure}

Our goal in this paper is to determine the homomorphism poset $\mathcal{K}_{2,n}$ of one family of complete bipartite graphs.  For small values of $n$, this is easy. Up to geo-isomorphim, there is only one realization of $K_{2,1}$ and so $\mathcal{K}_{2,1}$ is trivial.  There are only two realizations of $K_{2,2}$, one with no crossings and one with exactly one crossing.  The vertex labels in Figure \ref{fig:K22} indicate a geo-homomorphism that shows that $\mathcal{K}_{2,2}$ is a 2-element chain.
\medskip

For $n >2$, certainly a plane representation of $K_{2,n}$ will still be the first element of  $\mathcal{K}_{2,n}$, but the rest of the homomorphism poset is less obvious. To systematize our study, we develop a correspondence between geometric realizations of $K_{2,n}$ and permutations in $S_n$, defined in Section  2, in which edge crossings correspond to inversions.  
In Section 3, we give necessary and sufficient conditions for two permutations to correspond to geo-isomorphic realizations; we call such permutations {\it geo-equivalent}.  These conditions can be efficiently expressed using a directed version of permutation graphs.  The section includes a complete list of the geo-equivalence classes of $S_n$ for $n =4$ and $ 5$, as well as the number of geo-equivalence classes for all $n \leq 9$.
Some results on the size of geo-equivalence classes are given in Section 4, based on the structure of the modular decomposition tree of the permutation digraph. The poset structure of $\mathcal{K}_{2,n}$ is determined in Section 5, which includes Hasse diagrams for $\mathcal{K}_{2,4}$ and $\mathcal{K}_{2,5}$.  We compare the corresponding poset structure of the geo-equivalence classes of $S_n$ with that induced by the weak Bruhat order. We close with some open questions in Section 6.
\smallskip

Throughout this paper, the vertex set  of $K_{2,n}$ is denoted by $U=\{a,b\}$ and $V_n= \{1, 2, \dots, n\}$.  

\section{Permutations and Realizations of $K_{2,n}$}

\begin{figure}[h] 
   \centering
   \includegraphics[width=2.2in]{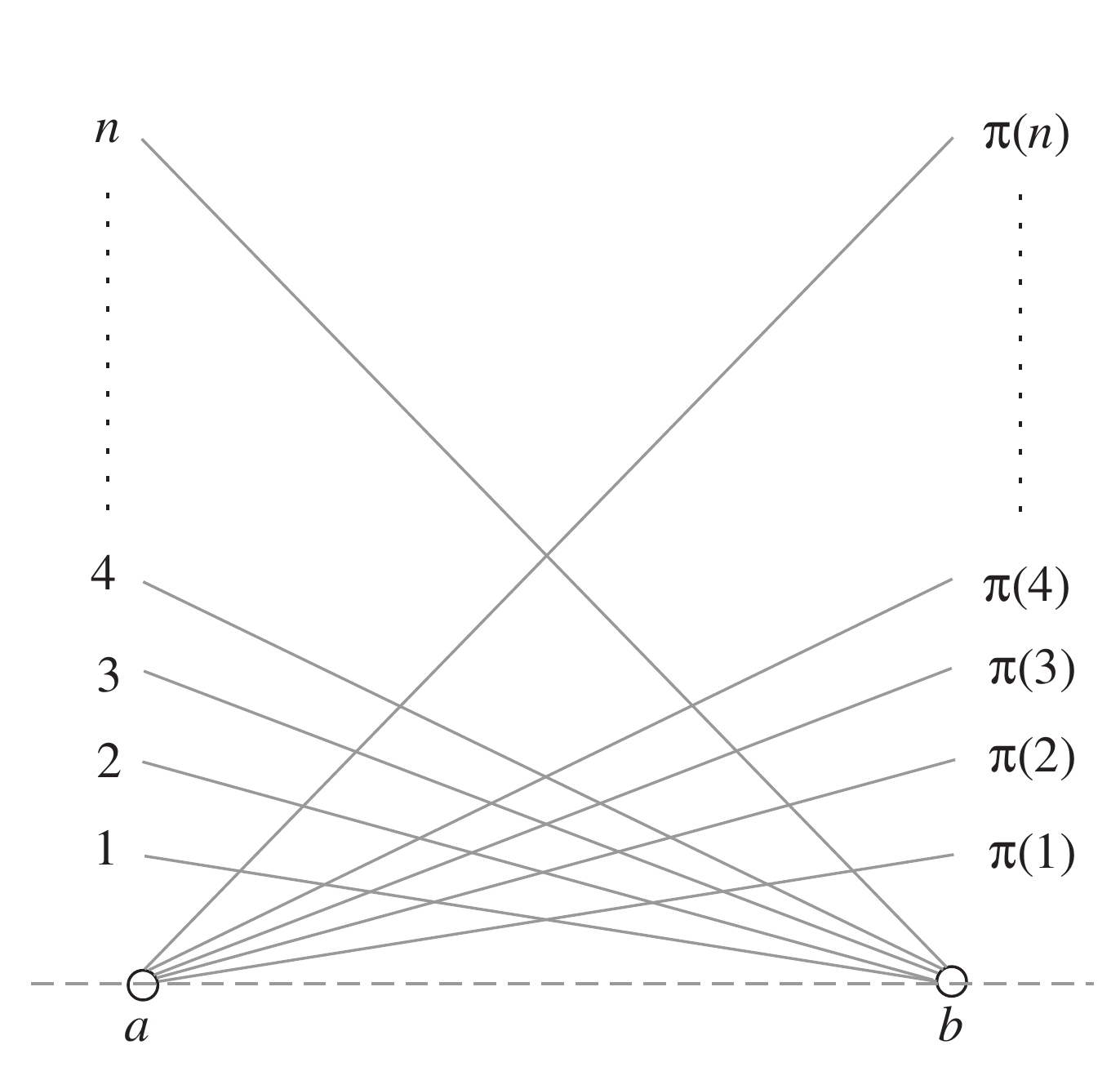} 
   \caption{Template for the construction.}
   \label{fig:template}
\end{figure}

For any  $\pi \in S_n$, we define a corresponding geometric realization of $\K2n$, denoted $\overline{K}_{2,n}(\pi)$, as follows.  We start with a template; from each of the points $a$ and $b$ in $\mathbb{R}^2$, draw $n$ intersecting rays, on the same side of the line $ab$. Label the rays emanating from $b$ consecutively $1$ through $n$; label the rays emanating from $a$ with $\pi(1)$ through $\pi(n)$, as in Figure \ref{fig:template}.
\medskip

\begin{figure}[h] 
   \centering
   \includegraphics[width=4.5in]{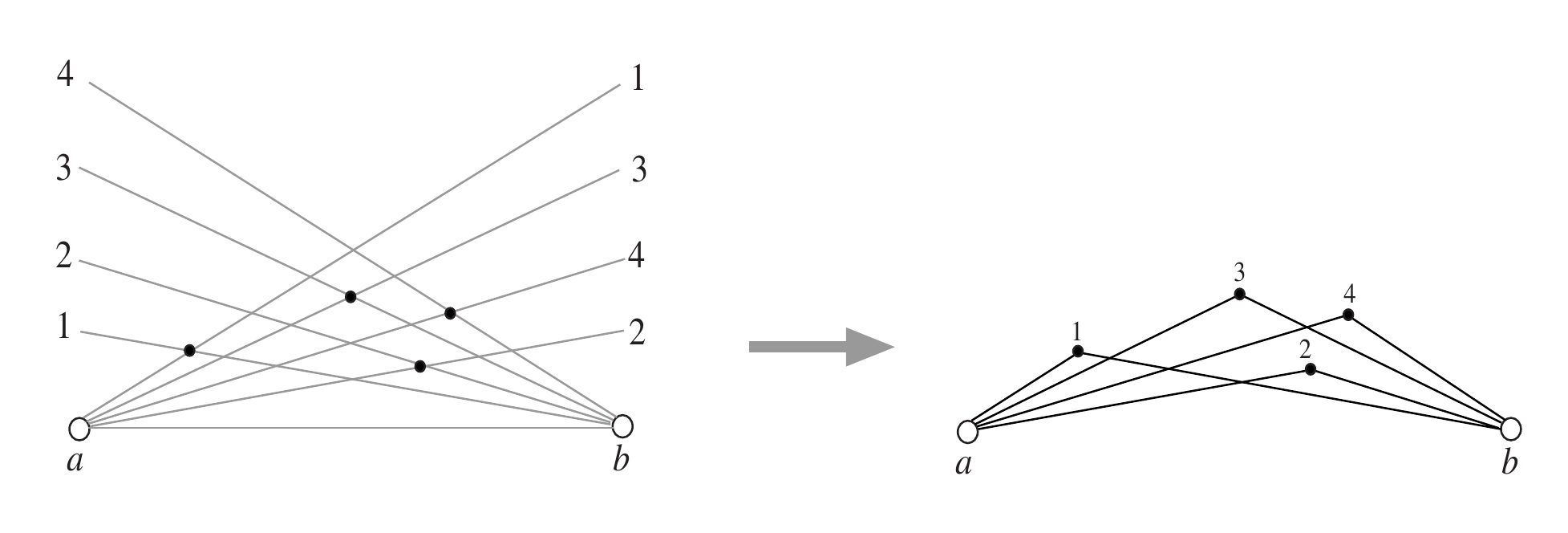} 
   \caption{The realization  $\overline{K}_{2,4}(2431)$.}
   \label{fig:easyeg}
\end{figure}

For each $i \in V_n$, position vertex $i$ at the intersection of the rays $ai$ and $bi$ on the template. With all the vertices in place, add the appropriate edges.  For example,  Figure \ref{fig:easyeg} illustrates the realization of $K_{2,4}$ corresponding to $\pi = 2431$. (We express permutations in word form, $\pi = \pi(1) \pi(2) \dots \pi(n)$, unless otherwise noted.)
\medskip

We would like to relate geometric properties of the realization $\Kpi$ to combinatorial properties of the permutation $\pi$.  To this end, recall that an inversion in a permutation is an instance of a smaller number appearing after a larger number.  For example, $2431$ contains exactly four inversions:  $1$ appears after $2, 3$ and $4$, and $3$ appears after $4$.  We state this definition more formally.

\begin{defn}
Let $i, j \in \{ 1, 2, \dots, n\}$ and let $\pi \in S_n$.  Then  $(i,j) $ is an {\it inversion} in $\pi$  if and only if $i < j $ and $\pi^{-1} (i) > \pi^{-1}(j)$. The set of inversions of  $\pi $ is denoted by $E(\pi)$ (also called the inversion set of $\pi$).
\end{defn}

A useful result that follows immediately from the definition  is
 \[
(i, j) \in E(\pi) \iff  \big(\pi^{-1}(j), \pi^{-1}(i)\big) \in E(\pi^{-1}),
\]
or equivalently, 
\[
(k,l) \in E(\pi^{-1}) \iff \big(\pi(l), \pi(k) \big) \in E(\pi).
\]  
\medskip

Returning to our example, we have 
$
E(2431) = \big\{(1, 2), (1, 3) , (1, 4), (3, 4)\big\}.
$
From Figure \ref{fig:easyeg}, we can see that in $\overline{K}_{2,4}(2431)$, $b1$ crosses $a2$, $a3$ and $a4$  and $b3$ crosses $a4$ . Moreover,  these are the only crossings.  This observation generalizes.

\begin{thm}\label{thm:inversion}
Let $\pi \in S_n$ and $i, j \in \{ 1, 2, \dots, n\}$. Then $bi$ crosses $aj$ in  $\Kpi$ if and only if $(i,j)\in E(\pi)$.
\end{thm}

\begin{figure}[htbp] 
   \centering
   \includegraphics[width=4.8in]{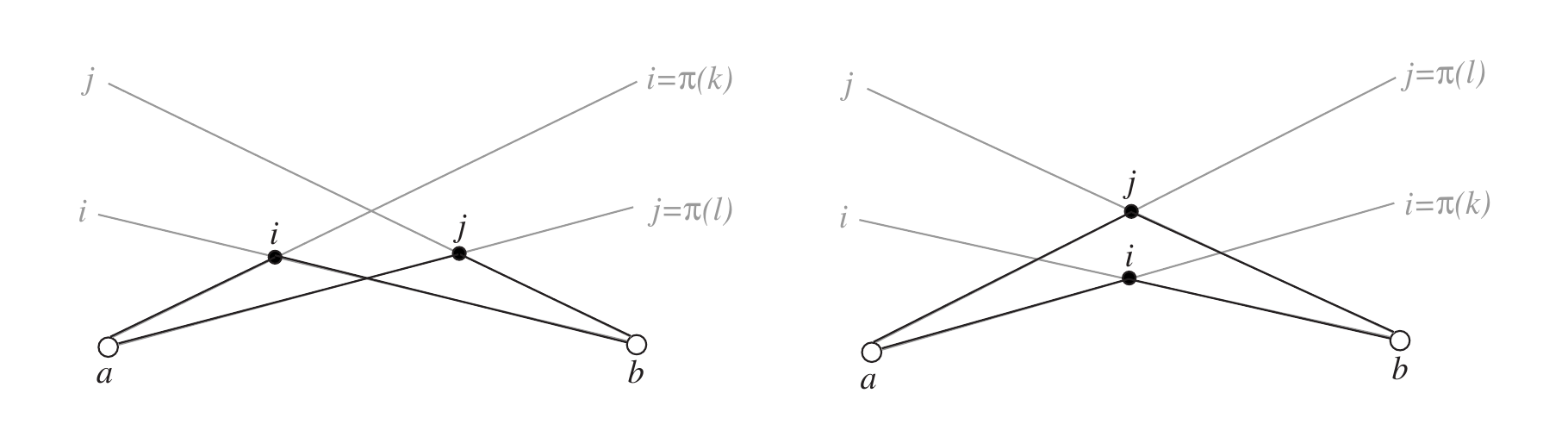} 
   \caption{Inversions correspond to crossings.}
   \label{fig:inversion}
\end{figure}

\begin{proof}
This result is obvious if we focus on the portion of the construction involving only vertices $i=\pi(k)$ and $j=\pi(l)$; see Figure \ref{fig:inversion}.  Up to geometric isomorphism, we get the subgraph on the left when  $k>l$, or equivalently, $\pi^{-1}(i) > \pi^{-1}(j)$, which by definition is when $(i,j)\in E(\pi)$. If $k < l$, then $\pi^{-1}(i) < \pi^{-1}(j)$, and so $(i, j) \not \in E(\pi)$.  In this case, we get the subgraph on the right. \end{proof}

Figure \ref{fig:inversion} also shows that if $i <j$, then $ai$ can never cross $bj$ in $\Kpi$ for any $\pi \in S_n$.

\begin{cor} \label{cor:total} Let $\pi, \sigma \in S_n$.
\begin{enumerate}
\item The total number of crossings in $\Kpi$ is $|E(\pi)|$.
\item If $\Ksig \cong \Kpi$, then $|E(\sigma)| = |E(\pi)|$.
\end{enumerate}
\end{cor}

\begin{proof}
The first statement follows immediately from Theorem \ref{thm:inversion}; the second follows from the fact that geo-isomorphisms preserve total number of crossings.
\end{proof}

\begin{figure}[h] 
   \centering
   \includegraphics[width=4in]{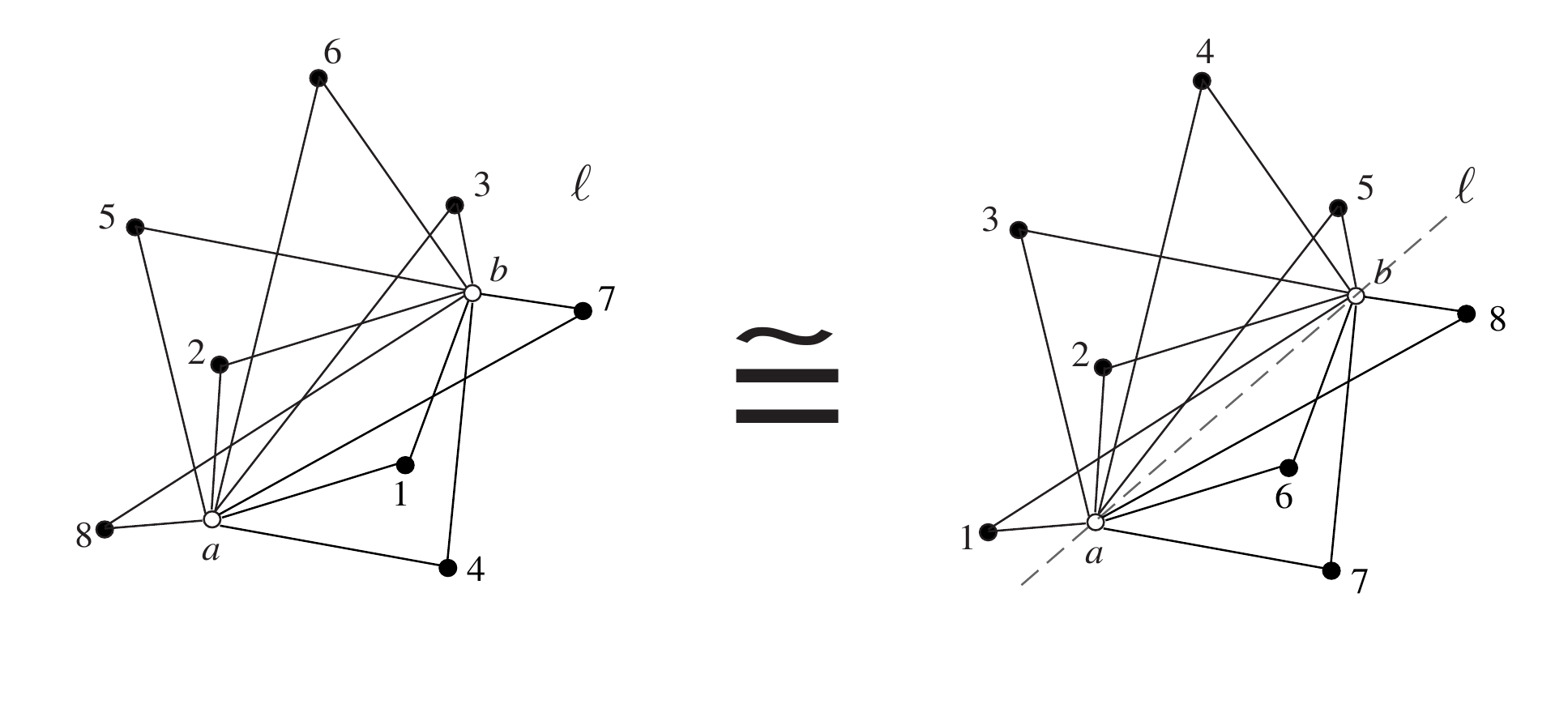} 
   \caption{Re-labeling vertices of a realization of $\overline{K}_{2,8}$.}
   \label{fig:gettoperm}
\end{figure}

Next, we show that any geometric realization of $\K2n$  is geo-isomorphic to $\Kpi$ for some $\pi \in S_n$.  Start with a (labeled) realization  $\overline{K}_{2,n}$ in the plane.  The line $\ell$ through points $a, b$ divides the plane into two half-planes.  Randomly select one half-plane and suppose  it contains $\{i_1, \dots i_t\} \subseteq V_n$ (where $0 < t \leq n$).  For each $1 \leq j \leq t$, let $\theta_b(j)$ denote the angle $\angle abi_j$.  Since the vertices are in general position, we can arrange these angles in strictly increasing order: 
\[
0 < \theta_b(j_1) < \theta_b (j_2) < \dots < \theta_b (j_t)  < 180^o.
\]
 Re-label vertex $i_{j_k}$ with $k$, so that  now $0 < \theta_b(1) < \theta_b (2) < \dots < \theta_b (t)  < 180^o$.  Next, let $\theta_a(j) = \angle baj$.  Arranging these angles in strictly increasing order induces a permutation of $\{1, \dots, t \}$, 
\[
0 < \theta_a\big (\pi(1)\big) < \theta_a \big(\pi(2)\big) < \dots < \theta_a \big(\pi(t)\big)  < 180^o.
\]
If $t = n$, we stop.  If $t < n$, 
re-label the remaining vertices $t+1, \dots, n$ so that $0 < \theta_b(t+1) < \theta_b (t+2) < \dots < \theta_b (n)  < 180^o$. Again, arranging the angles $\theta_a(j)$ in increasing order induces a permutation on $\{t+1, \dots, n \}$, 
\[
0 < \theta_a\big (\pi(t+1)\big) < \theta_a \big(\pi(t+2)\big) < \dots < \theta_a \big(\pi(n)\big)  < 180^o.
\]
Figure \ref{fig:gettoperm} illustrates the re-labeling protocol on a particular realization of $K_{2,8}$. The corresponding induced permutation is $\pi = 54231867$.

\begin{prop}\label{prop:Sn}
If the vertices of a geometric realization  $\overline{K}_{2,n}$ are re-labeled as above, with corresponding induced permutation $\pi$,  then  $\overline{K}_{2,n} \cong \Kpi$.
\end{prop}
  
  \begin{proof}
  By Theorem \ref{thm:inversion}, it suffices to show that for all $1 \leq i < j \leq n$,  $bi$ crosses $aj$  in  $\overline{K}_{2,n}$ if and only if $(i, j) \in E(\pi)$.\medskip
  
   First note that  if $i \leq t < j$, then by the re-labeling protocol, $i$ and $j$ are on opposite sides of line $\ell$, and so $bi$ cannot cross $aj$.  The construction forces $i = \pi(k)$ for some $k \in \{ 1, \dots, t \}$ and $j = \pi(l)$ for some $l \in \{t+1, \dots, n\}$.  Hence $k < l$, meaning $\pi^{-1}(i)  < \pi^{-1}(j)$, and so $(i,j) \not \in E(\pi)$.\medskip
  
 Next, assume $1 \leq i < j \leq t \, $ or $\,  t+1 \leq i < j $.  This means that $\theta_b(i) < \theta_b(j)$.  It is not difficult to see that  $bi$ crosses $aj$ if and only if $\theta_a(i) > \theta_a(j)$.  Letting $k=\pi^{-1}(i)$ and $l=\pi^{-1}(j)$, this is equivalent to $\theta_a\big(\pi(k)\big) > \theta_a\big(\pi(l) \big)$.  By construction, this occurs if and only if $k >l$.  By definition, this is true if and only if $(i,j) \in E(\pi)$.\medskip
 \end{proof}
  
\begin{figure}[h] 
     \centering
     \includegraphics[width=4.5in]{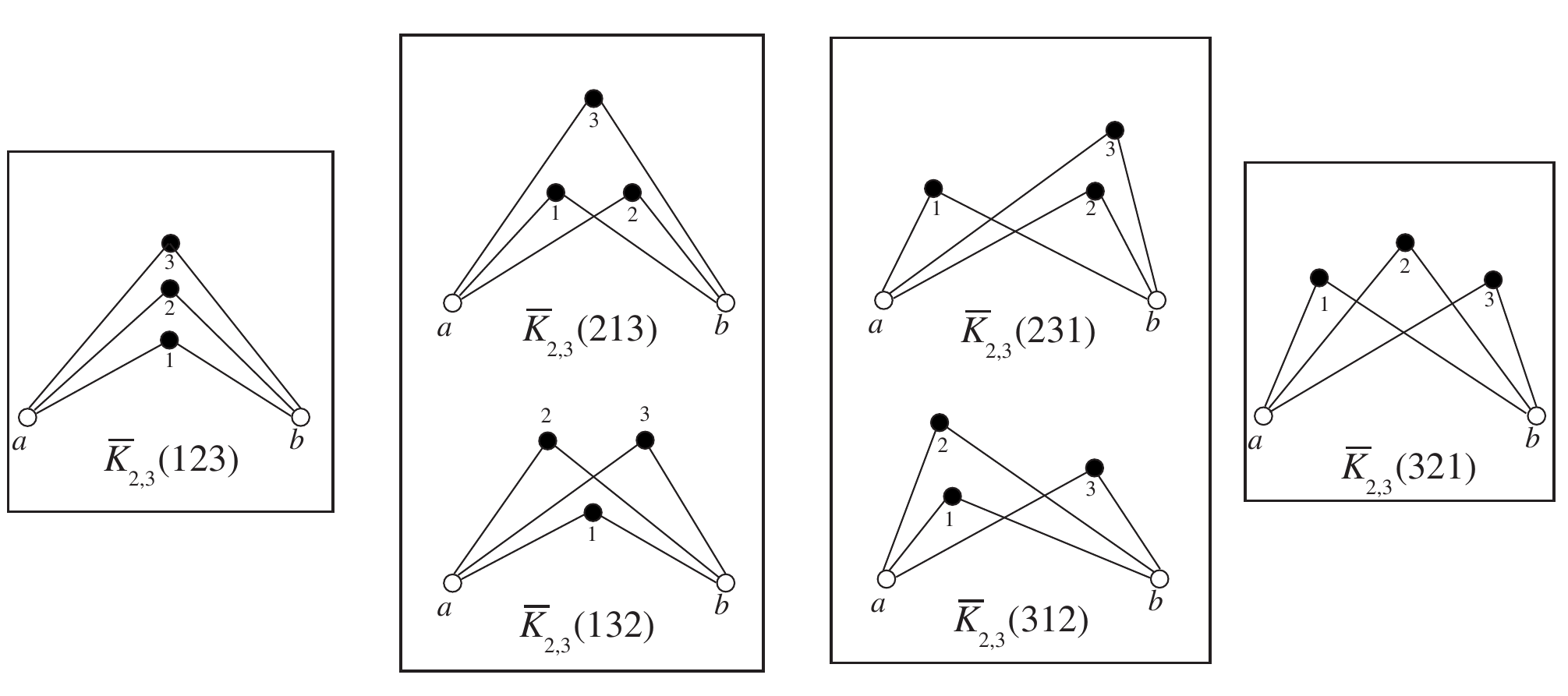} 
     \caption{Realizations of $K_{2,3}$.}
     \label{fig:K23elts}
  \end{figure}

  Applying this to the case $n = 3$, we conclude that the number of different geometric realizations of $K_{2,3}$ is at most $|S_3| = 6$. Visual inspection of  Figure~\ref{fig:K23elts} makes clear that in fact there are only 4 geo-isomorphism classes.  The question we address in the next section is: when do two permutations induce geo-isomorphic realizations?
  \medskip
  
Note that we can use geo-isomorphism to define an equivalence relation directly on $S_n$ by setting 
\[
\sigma \sim \pi \iff \Ksig \cong \Kpi.
\]
We will denote the {\it geo-equivalence} class of $\pi$ by $[\pi]$.  Similarly, we can define a partial order on the set of all geo-equivalence classes of $S_n$ by
\[
[\sigma] \preceq [\pi] \iff \Ksig \preceq \Kpi \text{ in } \mathcal{K}_{2,n}.
\]
We denote the resulting poset by $[\mathcal{S}_n]$.  From the above,  
$[\mathcal{S}_3]$ is the chain
\[
[123] \prec [213] \prec [231] \prec [321].
\]
\medskip

\section{Geo-equivalence Classes}\label{section:geo-classes}

In this section, we determine necessary and sufficient conditions for  two permutations to be geo-equivalent. We begin by defining an action of permutations on inversion sets.

\begin{defn}
Let $\sigma, \rho \in S_n$.  For all $(i,j) \in E(\sigma)$, let
\[
\rho \ast (i,j) = 
\begin{cases}
\big( \rho(i), \rho(j) \big) \quad &\text{ if } \rho(i) < \rho(j);\\
\big( \rho(j), \rho(i) \big)  &\text{ if } \rho(i) > \rho(j).
\end{cases}
\]
We say $\rho$ is {\it order-preserving} on $(i,j)$ in the first case, and  {\it order-reversing} on $(i,j)$ in the second.  We let $\rho \ast E(\sigma)$ denote the set $\{ \rho \ast (i,j) \, | \, (i,j) \in E(\sigma) \}$.
\end{defn}

The image of an inversion set under the action of a permutation may or may not itself be an inversion set.  We consider three illustrative examples.
 \medskip
 
\begin{ex} \label{ex:1}
If $\sigma_1 = 3214$ and $\rho_1 = 2341$, then 
 \begin{align*}
 E(\sigma_1)& =\{ (1,2), (1, 3), (2,3) \}\\
  \rho_1 \ast E(\sigma_1)& =\{ (2, 3), (2, 4), (3, 4) \} 
  = E(1432).
 \end{align*}
 Note that  $\rho_1$ is order-preserving on all inversions of $\sigma_1$.
 \end{ex}
 
 \begin{ex} \label{ex:2}  
 Let $\sigma_2 = 4312$ and $\rho_2 = 2341$; then 
 \begin{align*}
 E(\sigma_2) &= \{ (1, 3), (1,4), (2,3), (2, 4), (3, 4)\} \\
\rho_2 \ast E(\sigma_2) &= \{ (2, 4),  (3, 4), (1,4) \}  \cup \{ (1, 2), (1,3) \} \\
 & = E(4231).
 \end{align*}
In this case, $\rho_2$ is order-preserving on some inversions  of $\sigma_2$ and order-reversing on others.
\end{ex}

\begin{ex}\label{ex:3}   Let $\sigma_3 = 2413$ and $\rho_3 = 1324$; then 
 \begin{align*}
 E(\sigma_3) &= \{ (1, 2), (1, 4), (3, 4)\} \\
 \rho_3 \ast E(\sigma_3) &= \{ (1, 3), (1, 4),(2, 4) \}.
\end{align*} 
In this case, $\rho_3$ is order-preserving on all inversions of $\sigma_3$.  However,  the image is not the inversion set of any permutation.  To prove this, we need some additional background (an excellent overview of which can be found in Chapter 7 of  \cite{Gol80}).
\end{ex}

The inversions of a permutation $\pi \in S_n$ can be recorded in a graph $G(\pi)$, on vertices  $V_n = \{ 1, 2, \dots, n \}$, with  $i <j$ adjacent if and only if $(i,j) \in E(\pi)$.  More generally, we have the following definition.

\begin{defn}
A graph $G = (V, E)$ on $n$ vertices is a {\it permutation graph} if and only if there exists a bijection $L:V \to \{ 1, 2, \dots, n \}$ and a permutation $\pi \in S_n$ such that $L: G \to G(\pi)$ is a graph isomorphism.  In this case, we say $G$ {\it represents} $\pi$.
\end{defn}

Permutation graphs are related to another family of graphs, defined below.

\begin{defn}
A graph $G = (V, E)$ is {\it transitively orientable} if and only if  
its edges can be assigned an orientation $F$ so that in the directed graph $D=(V, F)$,   $(u,v),(v,w) \in F$ implies  $(u,w) \in F$.
\end{defn}
In 1971, Pneuli, Lempel and Even proved the following characterization of permutation graphs.

\begin{thm}\label{thm:ELP}\cite{ELP}
A graph $G$ is a permutation graph if and only if both $G$ and its complement $G^c$ are transitively orientable.
\end{thm}

We can rephrase this result in a way that allows us to quickly recognize when a set of ordered pairs is the inversion set of a permutation.

 \begin{cor}\label{cor:inversion_conditions}
Let $\mathcal{U}_n=\big\{(i,j) \mid \, 1 \leq i <j \leq n \big\}$, where $n \geq 2$, and let $A \subseteq \mathcal{U}_n$.  Then $A =E(\pi)$ for some $\pi \in S_n$ if and only if for all $i < j < k$,
\begin{enumerate}
\item $(i, j) \in A \text{ and }  (j,k) \in A \implies (i, k) \in A$;
\item $(i, j)  \in A^c \text{ and } (j, k) \in A^c \implies (i,k) \in A^c$.
\end{enumerate}
\end{cor}

An immediate consequence of this result is that the complement of an inversion set  in $\mathcal{U}_n$ is also an inversion set; 
in fact $[E(\pi)]^c = E(\pi^c)$, where $\pi^c$ is the `reverse' of $\pi$, given by
 $\pi^c = \pi(n) \pi(n-1) \dots \pi(1)$.
\medskip

Returning to Example \ref{ex:3}, we conclude that $\rho_3 \ast E(\sigma_3) = \{ (1, 3), (1, 4),(2, 4) \}$ is not an inversion set, as 
$
 (1,2), (2,3)   \in\big( \rho_3 \ast E(\sigma_3)\big)^c$, yet $(1, 3) \not \in \rho_3 \ast E(\sigma_3),
$ 
violating condition (2) of Corollary \ref{cor:inversion_conditions}.
Although the image of an inversion set is not always itself an inversion set, we do have the following result. 
\medskip

\begin{lem}\label{lem:rho_action}
For all $\sigma, \rho \in S_n$,  the image of $E(\sigma)$ under the action of $\rho$ is the symmetric difference, 
\[
\rho \ast  E(\sigma) = \big[E(\rho \cdot \sigma)\backslash E(\rho)\big] \cup \big[E(\rho) \backslash E(\rho \cdot \sigma)\big].
\]
More precisely,
\begin{align*}
E(\rho \cdot \sigma)\backslash E(\rho) &=\Big\{  \big( \rho(i), \rho(j) \big) \, | \, (i,j) \in E(\sigma) \text{ and } \rho(i) < \rho(j) \Big\}, \text{ and } \\
E(\rho) \backslash E(\rho \cdot \sigma) &= \Big\{\big( \rho(j), \rho(i) \big) \, | \, (i,j) \in E(\sigma) \text{ and }\rho(i) > \rho(j) \Big\}.
\end{align*} 
\end{lem}

\begin{proof}  
If $(k, l) \in E(\sigma)$, then $k < l$ and $\sigma^{-1}(k) > \sigma^{-1}(l)$.  
If $\rho(k) < \rho(l)$, then it is simply a matter of applying the definition to show that
$
\rho \ast (k,l) \in E(\rho\cdot \sigma) \backslash E(\rho).
$
Similarly, if $\rho(k) > \rho(l)$, then 
$
\rho \ast (k,l) \in E(\rho) \backslash E(\rho \cdot \sigma).
$
\medskip

Conversely, if $(i, j) \in E(\rho \cdot \sigma)\backslash E(\rho)$, then 
\[
\rho^{-1}(i) < \rho^{-1}(j) \text{ and } \sigma^{-1}\cdot \rho^{-1}(i) > \sigma^{-1} \cdot \rho^{-1}(j),
\]
 meaning that $(\rho^{-1}(i) , \rho^{-1}(j)) \in E(\sigma)$; clearly  
 $\rho \ast \big(\rho^{-1}(i) , \rho^{-1}(j)\big) = (i,j)$.
Similarly if $(i,j) \in E(\rho) \backslash E(\rho \cdot \sigma)$, then 
\[
\rho^{-1}(i) > \rho^{-1}(j) \text{ and }\sigma^{-1}\cdot \rho^{-1}(i) < \sigma^{-1} \cdot \rho^{-1}(j),
\]
 meaning that $(\rho^{-1}(j) , \rho^{-1}(i)) \in E(\sigma)$ and 
 $\rho \ast \big(\rho^{-1}(j) , \rho^{-1}(i) \big)= (i,j)$.
 \end{proof}

Note that if $\rho$ is order-preserving on all inversions of $\sigma$, then $\rho \ast E(\sigma)$ will never violate condition (1) of Corollary \ref{cor:inversion_conditions}.  For suppose $(k,l), (l,m) \in \rho \ast E(\sigma)$, where $k < l < m$.  Since $\rho$ preserves order, there exist $i<j<h$ such that 
\[
(k,l) = \big(\rho(i), \rho(j)\big), \, (l,m) = \big(\rho(j),  \rho(h)\big) \text{ and }(i,j), (j,h) \in E(\sigma).
\]  
Since $E(\sigma)$  satisfies (1), $(i, h) \in E(\sigma)$ and so  $\big(\rho(i),  \rho(h)\big) = (k,m) \in \rho \ast E(\sigma)$.  It can be shown similarly that  if $\rho$ is order-reversing on all inversions of $\sigma$, then $\rho \ast E(\sigma)$ satisfies condition (1).
\medskip
 
We are now ready for the main theorem of this section.
 
\begin{thm}\label{thm:complete}
Let $\sigma, \pi \in S_n$.  Then $\sigma \sim \pi$  if and only if there exists $\rho \in S_n$ such that
\begin{enumerate}
\item  $\rho \ast E(\sigma) = E(\pi)$;
\item $\rho$ is either order-preserving on $E(\sigma)$ or order-reversing on $E(\sigma)$.
\end{enumerate}
\end{thm}

\begin{proof}
Assume $\sigma \sim \pi$.  Then by definition, there exists a geo-isomorphism $f: \Ksig \to \Kpi$. 
Let  $\rho = f \mid_{V_n} \in S_n$.   First suppose $f(a)=a$ and $f(b) = b$.  
By Theorem \ref{thm:inversion} and the definition of geo-isomorphism, 
\begin{align*}
(i,j) \in E(\sigma) &\iff bi \text{ crosses }aj \text{ in } \Ksig \\
& \iff b\rho(i) \text{ crosses } a\rho(j) \text{ in } \Kpi \\
& \iff \big(\rho(i), \rho(j) \big) \in E(\pi).
\end{align*}
This implies both that $\rho$ is order-preserving on $E(\sigma)$ and that $\rho \ast E(\sigma) = E(\pi)$.
If $f(a)=b$ and $f(b) = a$, then $bi$ crosses $aj$ in $\Ksig$ if and only if $a\rho(i)$ crosses $b\rho(j)$ in $\Kpi$; in this case,  
 \[
 (i,j) \in E(\sigma) \iff \big (\rho(j), \rho(i) \big) \in E(\pi).
 \]
In this case,  $\rho \ast E(\sigma) = E(\pi)$ and $\rho$ is order-reversing on $E(\sigma)$.
\medskip

Conversely, assume  $\rho \ast E(\sigma) =  E(\pi)$ and  $\rho$ is order-preserving on $E(\sigma)$.
Define $g: \Ksig \to \Kpi$ by  
$
g(a) = a, \, g(b) = b$  and $g(i)=\rho(i)$  for all $ i \in V_n.
$
For $i <j$, we have
\begin{align*}
bi \text{ crosses } aj\text{  in } \Ksig & \iff (i, j) \in E(\sigma) \\
& \iff \big( \rho(i), \rho(j) \big) \in E(\pi)\\
& \iff b\rho(i) \text{ crosses } a\rho(j) \text{ in } \Kpi\\
& \iff g(b)g(i) \text{ crosses } g(a)g(j) \text{ in } \Kpi.
\end{align*}
Therefore $g$ is a geo-isomorphism.  If $\rho$ is order-reversing on $E(\sigma)$, then we  adapt this argument  by setting  
$
g(a) = b, \, g(b) = a \text{ and } g(i)=\rho(i) \text{ for all } i \in V_n.
$
\end{proof}

\begin{figure}[htbp] 
   \centering
   \includegraphics[width=4.5in]{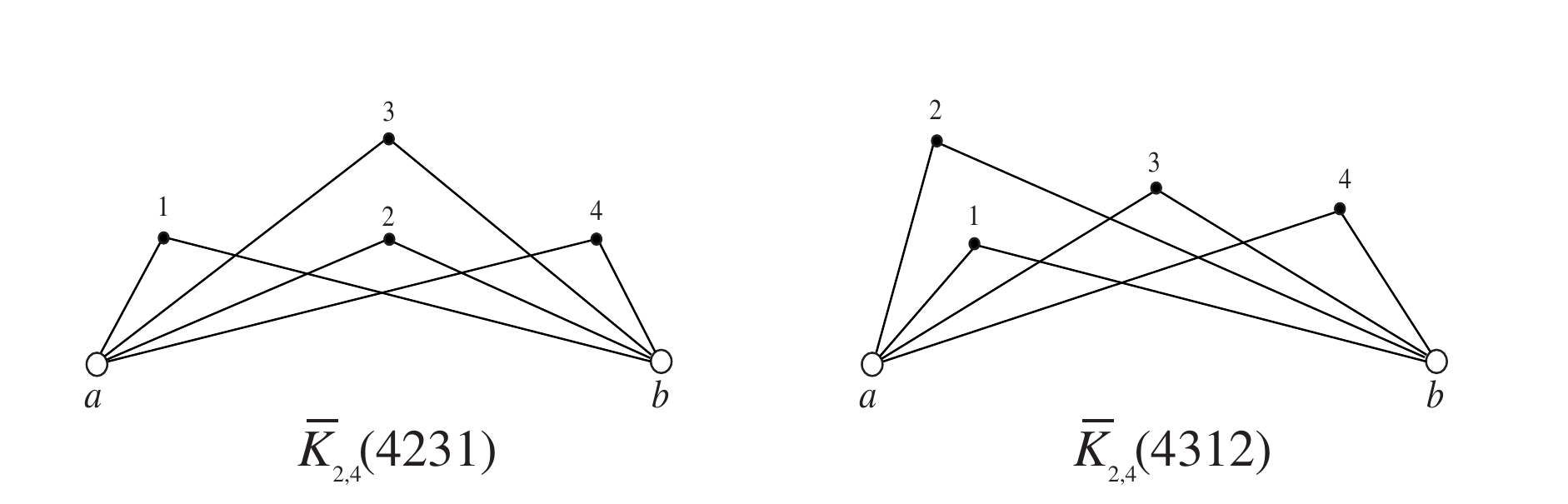} 
   \caption{Importance of condition (2) in Theorem \ref{thm:complete}.}
   \label{fig:countereg}
\end{figure}

Applying this theorem to Example~\ref{ex:1}, we conclude $1432 \sim 3214$, or equivalently, $\overline{K}_{2,4}(1432) \cong \overline{K}_{2,4}(3214)$.  Example~\ref{ex:2} illustrates the importance of condition (2) in Theorem \ref{thm:complete}; even though there exists  a $\rho \in S_n$ satisfying
 $
 \rho \ast E(4312) = E(4312),
$
the realizations of $K_{2,4}$ corresponding to  these two permutations (shown in Figure \ref{fig:countereg}) are not geo-isomorphic. One way to see this is to note that both edges incident to vertex 3 in $\overline{K}_{2,4}(4231)$ are crossed exactly once, but no vertex in $\overline{K}_{2,4}(4231)$ has this property.
\medskip

\begin{cor} \label{cor:inverse}
For all $\pi \in S_n$,  $\pi^{-1}\sim \pi$ via $\pi$, which is order-reversing.
\end{cor}

\begin{proof}  
This follows directly from our earlier observation that 
\[
(k,l) \in E(\pi^{-1}) \iff \big(\pi(l), \pi(k) \big) \in E(\pi).
\]
\end{proof}

\begin{ex}\label{ex:4}
Let $\pi = 3142$; then $\pi^{-1}=2413$.  In this case, $\pi^{-1} = \pi^c$, and so
\[
E(\pi^{-1}) = [E(\pi)]^c = \{ (1,2), (1, 4), (3,4) \}.
\]
We get $\pi \ast E(\pi^{-1}) = \{ (1,3), (2,3), (2,4)\} = E(\pi)$, with the action of $\pi$ reversing order on  
$E(\pi^{-1})$.  However, it is also true in this case that $\pi^{-1} \ast E(\pi^{-1}) = E(\pi)$, with the action of $\pi^{-1}$ preserving order on $E(\pi^{-1})$.
\end{ex}

\begin{cor} \label{cor:reverse}
For all $\pi \in S_n$, $\pi \sim ((\pi^c)^{-1})^c $ via $(\pi^c)^{-1}$, which is order-preserving.
\end{cor}

\begin{proof}
For $1 \leq i < j \leq n$, 
\begin{align*}
(i,j) \in E(\pi) & \iff (i,j) \not \in [E(\pi)]^c \iff (i,j) \not \in E(\pi^c) \\
&\iff 
(\pi^c)^{-1}(i) < (\pi^c)^{-1}(j).
\end{align*}
This shows that $(\pi^c)^{-1}$ is order-preserving on $E(\pi)$.
Next, for $1 \leq k < l \leq n$, 
\begin{align*}
(k,l) \in E \big (((\pi^c)^{-1})^c \big) & \iff (k,l) \not \in E((\pi^c)^{-1})\\
&\iff 
\pi^c (k) < \pi^c (l).
\end{align*}
Replacing $k$ with $(\pi^c)^{-1}(i)$ and $l$ with $(\pi^c)^{-1}(j)$, we get 

\begin{align*}
((\pi^c)^{-1}(i),(\pi^c)^{-1}(j)) \in E \big (((\pi^c)^{-1})^c \big) & \iff (\pi^c)^{-1}(i) < (\pi^c)^{-1}(j) \text{ and } i<j \\
& \iff (i,j) \in E(\pi).
\end{align*}
\end{proof}

 \noindent We can combine the last two corollaries to obtain the following.
 
\begin{cor}\label{cor:invrev}
For all $\pi \in S_n$, the permutations
$
\pi,\, \pi^{-1}, \,  ((\pi^c)^{-1})^c \text{  and  }   (((\pi^c)^{-1})^c)^{-1}
$
are all geo-equivalent.
\end{cor}

\begin{ex} By Corollary~\ref{cor:invrev},
\[
\pi = 2431, \, \pi^{-1}=3142, \,  ((\pi^c)^{-1})^c =3241\, \text{  and  }   ((\pi^c)^{-1})^c)^{-1}=4213
\]
are all geo-equivalent. 
\medskip

The four permutations in Corollary \ref{cor:invrev} may not all be distinct.  In Example~\ref{ex:4}, we saw that for $\pi =3142$, $\pi^{-1}= 2413 = \pi^c$,  so 
\[
((\pi^c)^{-1})^c = \pi^{-1} \text{  and  } (((\pi^c)^{-1})^c)^{-1}= \pi.
\]
For $\pi = 3412$, we have 
$
\pi = \pi^{-1} =  ((\pi^c)^{-1})^c =  (((\pi^c)^{-1})^c)^{-1}.
$
\end{ex}

Recall that for all $\pi \in S_n$,   $G(\pi)$ has vertices $V_n = \{1, 2, \dots, n\}$ and edges $\{ ij \, | \, (i,j) \in  E(\pi)\}$.  Thus $G(\sigma)$ is isomorphic to $G(\pi)$ if and only if there exists $\rho \in S_n$ such that $\rho \ast E(\sigma) = E(\pi)$; that is, $G(\sigma) \cong G(\pi)$ as abstract graphs if and only if condition (1), but not necessarily condition (2),  of Theorem~\ref{thm:complete} is satisfied.  In particular, Example 3 shows that $G(4312) \cong G(4231)$, yet $4312 \not \sim 4231$.  
Thus the number of geo-equivalence classes of $S_n$ may exceed the number of non-isomorphic permutation graphs on $n$ vertices.
\medskip

We can rephrase Theorem~\ref{thm:complete} in the language of permutation graphs by introducing 
a directed version of $G(\pi)$.
More precisely, 
for all $\pi \in S_n$,  we let $D(\pi)$ denote the digraph with vertex set  $V_n = \{1, 2, \dots, n\}$ and arc set $E(\pi)$.  We will call a digraph $D=(V,F)$ a {\it permutation digraph} if and only if $D \cong D(\pi)$ for some permutation $\pi$; in this case, we say $D$ {\it represents} $\pi$.

\begin{lem}\label{lem:invdigraph}
 Let $D = (V,F)$ be a digraph and let $-D = (V, -F)$ denote the digraph obtained by reversing direction on all arcs of $D$.  If $D \cong D(\pi)$ for some $\pi \in S_n$, then $-D \cong D(\pi^{-1})$.
\end{lem}

\begin{proof}  
Suppose $L : V \to \{1, 2, \dots, n \}$ is a bijection establishing $D \cong D(\pi)$; that is,
\[
(u,v) \in F \iff \big( L(u), L(v) \big) \in E(\pi).
\]
By Corollary  \ref{cor:reverse}, $(i,j) \in E(\pi) \iff (\pi^{-1}(j), \pi^{-1}(i) \big) \in E(\pi^{-1})$. Hence  the bijection $\pi^{-1} \circ L:V \to \{1, 2, \dots, n\}$  establishes $ -D \cong D(\pi^{-1})$.
\end{proof}

\begin{thm}\label{thm:permdigraphs}
Let $\pi, \sigma \in S_n$.  Then $\sigma \sim \pi$  if and only if  
either $D(\sigma) \cong D(\pi)$ or $D(\sigma) \cong D(\pi^{-1})$.
\end{thm}

\begin{proof}
If there exists $\rho \in S_n$ such that $\rho \ast E(\sigma) = E(\pi)$, with $\rho$ preserving order on $E(\sigma)$, then $\rho$ is also a digraph isomorphism $D(\sigma)\to D(\pi)$.  If $\rho$ is order-reversing 
on $E(\sigma)$, then $\pi^{-1} \circ \rho: D(\sigma) \to D(\pi^{-1}) $ is a digraph isomorphism.
Conversely,  a digraph isomorphism $\gamma:D(\sigma) \to D(\pi)$ must be an element of $S_n$ satisfying $\gamma \ast E(\sigma) = E(\pi)$, with $\gamma$ preserving order on $E(\sigma)$. If the digraph isomorphism is $\gamma:D(\sigma) \to D(\pi^{-1})$, then $(\pi \circ \gamma) \ast E(\sigma) = E(\pi)$, with $\pi \circ \gamma$ reversing order on $E(\sigma)$.
\end{proof}

\begin{figure}[htbp] 
   \centering
   \includegraphics[width=4.5in]{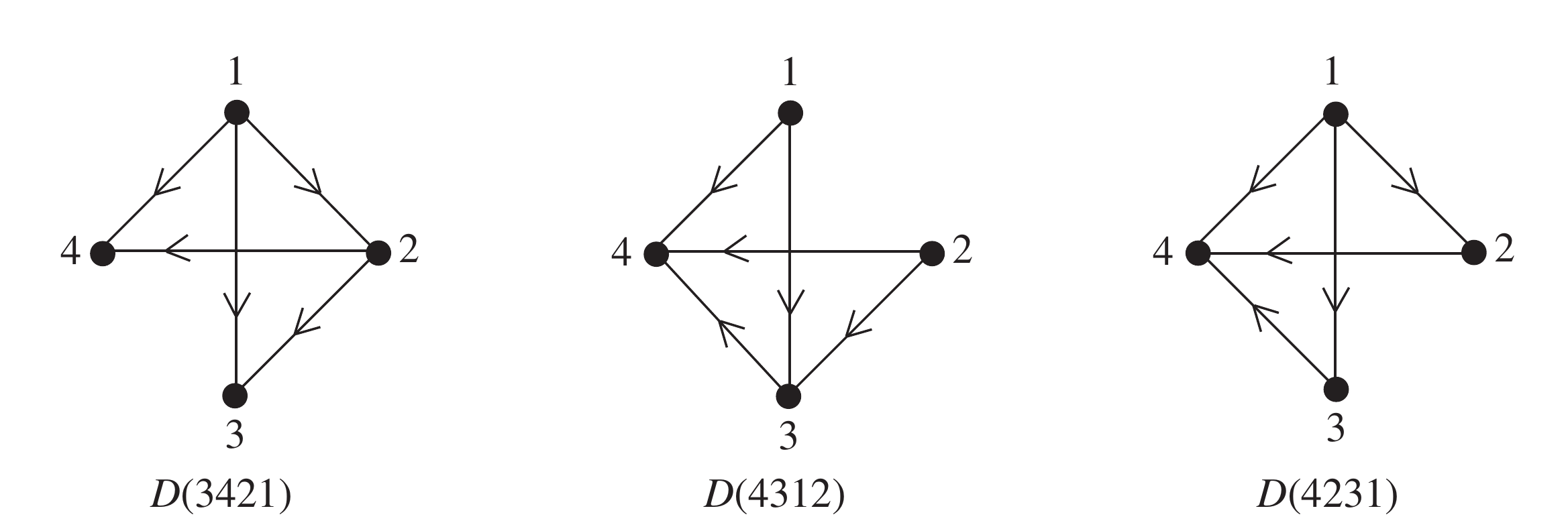} 
   \caption{Some permutation digraphs.}
   \label{fig:2transor}
\end{figure}

 Figure \ref{fig:2transor} illustrates the previous two results.  First note that the underlying undirected graphs are the same, so $G(3421) \cong G(4312) \cong G(4231)$; up to isomorphism, there is only one permutation graph on 4 vertices with 5 edges. 
 The first two digraphs, $D(3421)$ and $D(4312)$, are `reverses' of each other, as expected from the fact that $(3421)^{-1} = (4312)$.  The third digraph, $D(4231)$, is not isomorphic to either of the previous two.  Note that reversing the direction on all arcs of $D(4231)$ yields a digraph  isomorphic to the original, as expected from the fact that $(4231)^{-1} = 4231$. We conclude that the permutations of $S_4$ with $5$ inversions divide into two geo-equivalence classes: $[3421] = \{3421, 4312\}$ and $[4231] = \{ 4231\}$.
 \medskip
 
 Theorem \ref{thm:permdigraphs} suggests that  to determine the geo-equivalence classes of $S_n$, we must determine the isomorphism classes of permutation digraphs on $n$ vertices, and additionally identify a digraph $D$ with its reverse, $-D$.  Following Colbourn (see \cite{Col}), we call a permutation graph {\it uniquely orientable} if and only if it admits only one transitive orientation and its reverse.  Furthermore, we make the following definition.

\begin{defn}
Two permutations digraphs $D_1$ and $D_2$  are {\it related} if and only if either $D_1 \cong D_2$ or $D_1 \cong -D_2$; otherwise they are unrelated.
\end{defn}

Using this terminology, the number of geo-equivalence classes of $S_n$ is the number of unrelated permutation digraphs on $n$ vertices.  Table~\ref{table:S4}  gives the partitioning of $S_4$ into geo-equivalence classes.  There are $11$ non-isomorphic (undirected) graphs on $4$ vertices and all of them are permutation graphs.  The only one that is not uniquely orientable is the one in Figure \ref{fig:2transor}, giving 12 geo-equivalence classes in total.
\medskip

\begin{table}[h]
\caption{Geo-equivalence classes of $S_4$.}
\begin{center}
\begin{tabular}{c}
\includegraphics[scale=0.55]{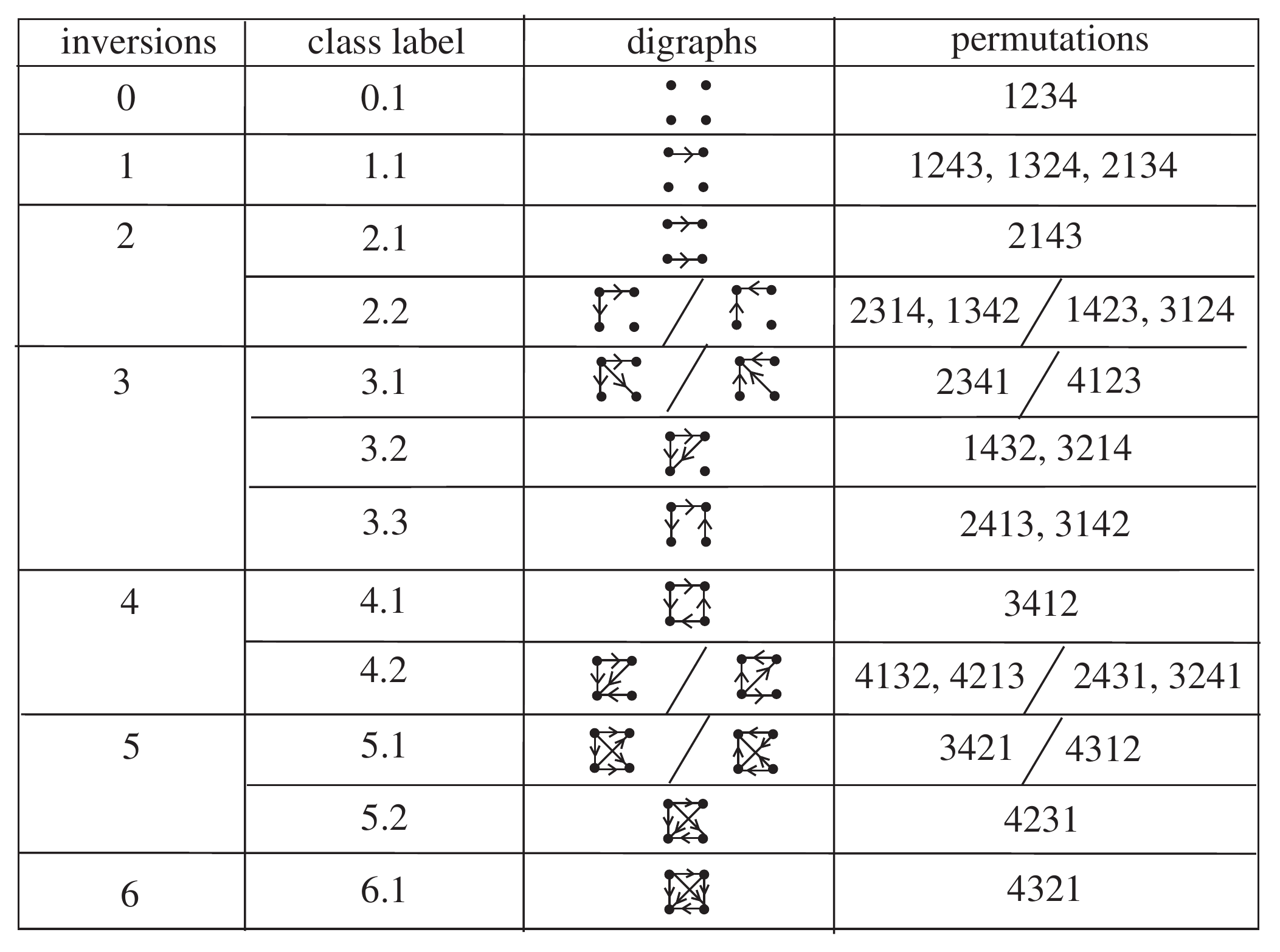}
\end{tabular}
\end{center}
\label{table:S4}
\end{table}%

\begin{figure}[h] 
   \centering
   \includegraphics[width=4.8in]{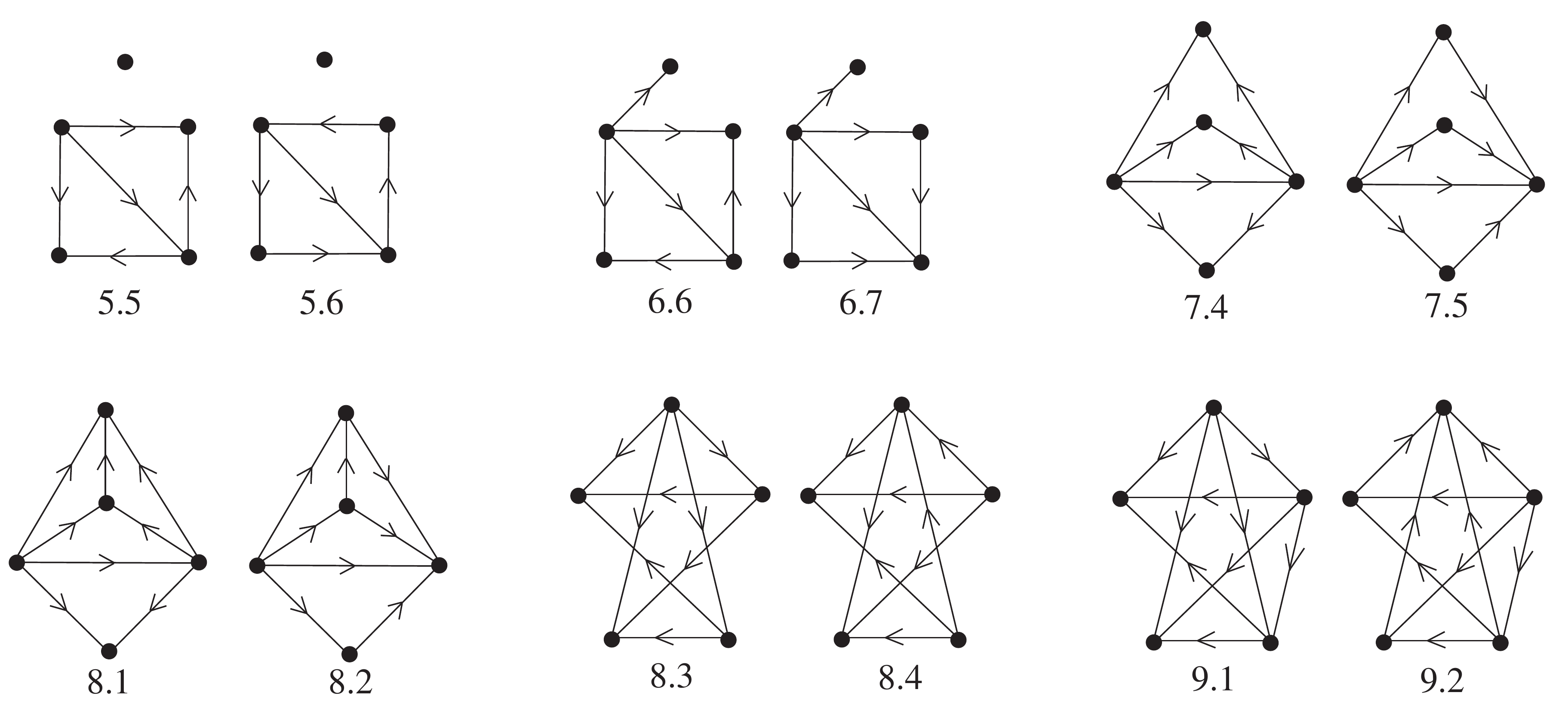} 
   \caption{Six  permutations graphs, $n=5$, each with 2 unrelated transitive orientations.}
   \label{fig:unrelatedn5}
\end{figure}

Progressing to $n=5$, there are $34$ non-isomorphic graphs in total, but one of them,  $C_5$, is not transitively orientable and is therefore not a permutation graph.  Of the remaining $33$ graphs, $27$  are uniquely orientable and the remaining $6$ have exactly two unrelated orientations, as shown in Figure~\ref{fig:unrelatedn5}.  Thus, $S_5$ has $39$ geo-equivalence classes in total; these are given in Table~\ref{table:S5}. (To save space, this table does not include diagrams of each possible oriented digraph.) As in the previous table, permutations corresponding to opposite orientations of the underlying graph are separated by a diagonal slash.
\medskip

\medskip

\begin{center}
\renewcommand{\arraystretch}{1.5}
\tablefirsthead{%
	\hline
	 inversions	   & class label 	     & permutations				 \\ 
	\hline}
\tablehead{%
	\hline
	\multicolumn{3}{|l|}{\small\sl continued from previous page}\\
	\hline
	 inversions	& class label 	& permutations 				 \\
	 \hline}
\tabletail{%
	\hline
	\multicolumn{3}{|r|}{\small\sl continued on next page}\\	
	\hline}
\tablelasttail{\hline}
\tablecaption{Geo-equivalence classes of $S_5$.}\label{table:S5}
\begin{supertabular}{|c|c|c|}
0					&0.1		& 12345							\\ \hline
1					&1.1		& 12354, 12435, 13245, 21345		\\ 	 \hline

2					&2.1		& 12453, 13425, 23145 / 	12534, 14235, 31245		\\   \cline{2-3}
					&2.2		&13254, 21354, 21435				\\ \hline

3					&3.1		& 13452, 23415 /  15234, 	41235		\\ \cline{2-3}
					&3.2		& 13524, 24135 / 14253, 31425		\\ \cline{2-3}
					&3.3		&12543, 14325, 32145							\\ \cline{2-3}
					&3.4		&21453, 23154 / 21534, 31254	\\ \hline

4					&4.1		& 23451	/ 51234	\\ \cline{2-3}
					&4.2		& 13542, 14352, 24315, 32415 /  15234, 15324, 41325, 42135 \\ \cline{2-3}
					&4.3		& 23514, 31452 / 41253, 25134		\\  \cline{2-3}
					&4.4		& 24153  / 31524	\\ \cline{2-3}
					&4.5		&14523, 	34125		\\ \cline{2-3}
					&4.6		&32154, 	21543					\\ \hline

5					&5.1		& 23541, 24351, 32451 /	51243, 51324,  52134		\\  \cline{2-3}

					&5.2		& 34152, 24513 /  35124, 41523 		\\  \cline{2-3}
					&5.3		& 25314 / 41352	\\ \cline{2-3}
					&5.4		& 32514, 31542 /  42153, 	25143	\\ \cline{2-3}
					&5.5		& 14532, 34215 / 15423, 43125	\\ \cline{2-3}
					&5.6		&15342, 	42315				\\ \hline

6					&6.1		& 32541 / 52143	\\ \cline{2-3}
					&6.2		& 34512 / 45123	\\ \cline{2-3}
					&6.3		& 25413, 43152 / 41532, 35214		\\ \cline{2-3}
					&6.4		&15432, 	43215							\\ \cline{2-3}
					&6.5		&35142, 	42513				\\ \cline{2-3}
					&6.6		& 24531, 34251 / 51423, 53124 	\\  \cline{2-3}
					&6.7		& 25341, 42351 / 51342, 52314 	\\ \hline

7					&7.1		& 25431, 43251  /  51432, 53214		\\  \cline{2-3}
					&7.2		& 35412, 43512  / 45132, 45213		\\ \cline{2-3}
					&7.3		& 42531,  35241 / 52413, 	53142	\\ \cline{2-3}
					&7.4		& 34521	/ 54123 	\\ \cline{2-3}
					&7.5		&52341							\\ \hline		

8					&8.1		& 35421, 43521 / 54132, 54213		\\ \cline{2-3}
					&8.2		&52431, 	53241					\\ \cline{2-3}
					&8.3		&45231 / 53412	\\  \cline{2-3}
					&8.4		&45312								\\ \hline

9					&9.1		& 45321/ 54312 	\\ \cline{2-3}
					&9.2		& 53421 / 54231	\\ \hline

10					&10.1		&54321								\\ \hline
\end{supertabular}
\end{center}
 
From the geo-equivalence classes for $n=3, 4$ and $5$, one might conjecture that involutions (\emph{i.e.} permutations of order two) can only be geo-equivalent to other involutions. However, a counterexample exists at $n=6$;  $\pi = 465132$ is an involution (in cycle notation, $\pi = (14)(26)(35)$), $\sigma = 465213$ is  the 6-cycle $(142635)$, yet 
$\sigma \sim \pi$ via $\rho = 231456$.
\medskip

We can write a program based on Theorem~\ref{thm:complete} to determine geo-equivalence classes for larger values of $n$.  If we let $a_n$ denote the number of geo-equivalence classes in $S_n$ (where $n \geq 1$), then the first nine terms of the integer sequence $(a_n)$ are:
\[
1, \, 2, \, 4, \,12, \,39, \, 182, \, 1033 , \, 7605, \, 66302,...
\]
Interestingly, this does not match any other sequence in the Online Encyclopedia of Integer Sequences.  However, implementing this theorem involves testing $n!$ permutations as candidates for $\rho$, and it is therefore very inefficient.  (However, the interested reader may find both  C++ and Python code for this algorithm at entry A180487 in OEIS \cite{oeis}.)
\medskip

For an approach based on Theorem~\ref{thm:permdigraphs}, we can start with $p_n \leq a_n$, where $p_n$ is the number of permutation graphs on $n$ vertices.  However, neither a closed nor a recursive formula for $p_n$ is known.   Evens, Lempel and Pnueli \cite{ELP} gave a polynomial-time algorithm for recognizing permutation graphs in 1971, and ten years later, Colbourne \cite{Col} gave a polynomial-time algorithm for determining if two permutation graphs are isomorphic.
More recently,  progress has been made on the enumeration of certain subclasses of permutation graphs.  
In 1999, Guruswami \cite{Guru} gave a generating function for the number of non-isomorphic cographs and threshold graphs.  (We will discuss cographs further in the next section.) Koh and Ree \cite{KohRee} found a recurrence relation for the number of vertex-labeld connected permutation graphs in 2007, and 
in 2009, Saitoh, Otachi, Yamanaka and  Uehara \cite{SOYU} developed a linear time algorithm for generating and enumerating  non-isomorphic bipartite permutation graphs. In the absence of a starting point for the number of geo-equivalence classes, we turn to determining the size of geo-equivalence classes.

\section{Size of Geo-equivalence Classes}

In this section, we develop a method for determining the size of the geo-equivalence class  represented by a given permutation digraph.
A useful tool for this investigation is modular decomposition, which we briefly review below.  Although this theory can be traced  back to a seminal 1967 paper by  Gallai  \cite{Gallai}, we use the more modern terminology and notation that can be found in Brandstadt, Le and Spinrad \cite{BLS} or McConnell \cite{McC}.
 
 \begin{defn}
 \hfill
 \begin{enumerate}
 \item A {\it module} of a graph $G=(V,E)$ is a set of vertices $M$ such that any vertex outside $M$ is either adjacent to every vertex in $M$, or not adjacent to any vertex in $M$.  More formally, for all $v \in V \backslash M$, either $uv \in E$ for all $u \in M$, or $uv \in E^c$ for all $u \in M$.
 \item Two modules $M$ and $N$ {\it overlap} if and only if $M \cap N$, $M \backslash N$ and $N \backslash M$ are all non-empty.
 \item A module $M$ is {\it strong} if and only if it does not overlap with any other module of $G$; otherwise it is {\it weak}.
 \end{enumerate}
 \end{defn}
 
For any graph $G$,  $V$ and $\{v\}$ for all $v \in V$ are modules (in fact, strong modules); they are called {\it trivial} modules.  Note that the modules (and strong modules) of $G$ and $G^c$ are the same.  
 \medskip

We can recursively partition the vertex set of a graph $G = (V,E)$ into its strong modules using the following algorithm.  We use $G|M$ to denote the subgraph of $G$ induced by $M$.  We begin the algorithm with $M = V$. 
\begin{enumerate}
\item If $|M| = 1$, then stop.
\item If $G|M$ is disconnected, then partition $M$ into its connected components.
\item If  $G|M$ is connected, but $G^c|M$ is disconnected, partition $M$ into the connected components of $G^c|M$.
\item If both $G|M$ and $G^c|M$ are connected, then $M$ can be partitioned into its maximal submodules (which will be strong modules of $G$).
\end{enumerate} 

The {\it modular decomposition tree} of $G$ has the strong modules of $G$ as its nodes, with $V$ being the root node, and the children of a node $M$ being the strong modules in the partition of $M$ from the algorithm above.  Every leaf in this tree is a singleton set, $\{v\}$. An internal node $M$ in the tree is called:
\begin{itemize}
\item a {\it degenerate 0-node} if  $G|M$ is disconnnected;
\item a {\it degenerate 1-node} if   $G|M$ is connected, but $G^c|M$ is disconnected;
\item a {\it prime} node if both $G|M$ and $G^c|M$ are connected.
\end{itemize}
Every weak module of $G$ is a union of children of degenerate node, and conversely, every union of children of a degenerate node is a weak module.   Note that a degenerate 0-node of $G$ is a degenerate 1-node of $G^c$, and vice versa.
\medskip

To every internal node $M$ of the modular decomposition tree of $G=(V,E)$, we associate a {\it quotient graph} $Q(M)$, whose vertices are the children of $M$, with two children $X$ and $Y$ being {\it adjacent} if and only if $xy \in E$ for some $x \in X$ and $y \in Y$.  Note that by definition of a module,  $xy \in E$ for some $x \in X, y \in Y$ if and only if $xy \in E$ for all $x \in X, y \in Y$.  It follows directly from the definitions that if $M$ is a degenerate 0-node, then $Q(M)$ is a null graph,
and if $M$ is a degenerate 1-node with $k$ children, then $Q(M)$ is a complete graph on $k$ vertices.
\medskip

\begin{lem}\label{lem:childorient} \cite{Gallai} \cite{McC}
If $X$ and $Y$ are adjacent children of either a degenerate 1-node or a prime node, then in any transitive orientation $F$ of $G$, all edges between $X$ and $Y$ must be oriented the same way ({\it i.e.} either $X \times Y \subseteq F$ or $Y \times X \subseteq F$).
\end{lem}

Thus any transitive orientation on the edges of $G$ unambiguously restricts to a transitive orientation on each quotient graph $Q(M)$.  Conversely, Gallai showed that any set of transitive orientations on the quotient graphs extends in the obvious way to a transitive orientation on $G$.  Since complete graphs are always transitively orientable, we conclude that $G$ is transitively orientable if and only if for every prime node $M$ of $G$, $Q(M)$ is transitively orientable.

\begin{prop}\cite{Gallai}\label{prop:prime}
Let $M$ be a prime node of the modular decomposition tree of $G$.  If $Q(M)$ is transitively orientable, then it is uniquely orientable.
\end{prop}

Putting all of these facts together, we can determine the number of different transitive orientations on a vertex-labeled transitively orientable graph.

\begin{cor}\label{cor:countTO}
Let $G=(V,E)$ be a transitively orientable graph.  Suppose the internal nodes of the modular decomposition tree of $G$ consist of:
\begin{itemize}
\item prime nodes $P_1, P_2, \dots , P_s$;
\item degenerate 1-nodes $N_1, \dots, N_t$, where $N_i$ has $k_i$ children;
\item degenerate 0-nodes $M_1, \dots, M_r$.
\end{itemize}
Then $G$ has $2^s \cdot k_1! \cdot k_2! \cdots k_t!$ different transitive orientations.
\end{cor}

Note that this number counts any transitive orientation $F$ and its reverse $-F$ as different orientations: isomorphic orientations are also counted as different. Hence, this is not the number of unrelated transitive orientations on $G$, only an upper bound.
\medskip

Recall that a graph $G$ is a permutation graph if and only if both $G$ and $G^c$ are transitively orientable.  In \cite{ELP}, Evens, Lempel and Pnueli give an algorithm that takes as input transitive orientations $F, F_1$ on  $G, G^c$ respectively, and outputs a permutation $\pi$ such that $(V, F) \cong D(\pi)$ ({\it i.e.} a permutation represented by $(V,F)$). 
First they show that superimposing the two orientations yields a transitively oriented complete graph, $(V, F \cup F_1)$.  
Associated with this orientation on $K_n$ is a unique vertex labeling function $L:V \to \{ 1, \dots, n \}$ satisfying
\begin{equation}\label{eq:defL}
L(v) < L(w) \iff (v, w) \in F \cup F_1.
\end{equation}
Since $F$ and $F_1$ are both transitive, by Corollary  \ref{cor:inversion_conditions} there exists a unique permutation $\pi$ such that 
\begin{equation}\label{eq:defpi}
\Big\{\big(L(v), L (w)\big)  \mid (v,w) \in F\Big \} = E(\pi).
\end{equation} 
We say $\pi$ is the permutation {\it induced by} $F \cup F_1$, or  equivalently, $F\cup F_1$ {\it induces} $\pi$.   
Note that we also have
\begin{equation}\label{eq:defpic}
\Big\{\big(L(v), L (w)\big)  \mid (v,w) \in F_1 \Big \} = E(\pi)^c = E(\pi^c).
\end{equation} 
Hence we have both $L:(V, F) \cong D(\pi)$ and $L:(V, F_1) \cong D(\pi^c)$.
\medskip

Given a permutation digraph $D=(V, F)$, with underlying undirected graph $G$, this algorithm defines a function 
\[
\Phi:\{ \text{transitive orientations on } G^c\} \to \{ \text{permutations represented by } D\}.
\]
Now $\Phi$ is surjective, for assume  $\pi$ is represented by $D= (V, F)$. By definition, there exists a bijection $L: V \to \{ 1, \dots, n \}$ that is 
an isomorphism $D \to D(\pi)$.  Applying $L^{-1}$ to the vertices of $D(\pi^c)$ 
induces a transitive orientation $F_1$ on $G^c$, and it is clear that $F \cup F_1$ will induce $\pi$. 
However, 
$\Phi$ is not injective, as can be seen by letting $D$ be the null digraph on $n$ (labeled) vertices. In this case, $G^c$ is complete, and has $n!$ different transitive orientations. However, the only permutation $D$ represents is the identity in $S_n$.    The following proposition gives a necessary and sufficient condition for $\Phi$ to  take two transitive orientations of $G^c$ to the same permutation.

\begin{prop}\label{prop:Phinotinjective}
Let $D = (V, F)$ be a permutation digraph with underlying undirected graph $G$, and let $F_1, F_2$ be two transitive orientations on $G^c$.  Then $F \cup F_1$ and $F \cup F_2$ induce the same permutation if and only if there exists a bijection $f:V \to V$ such that $f:(V, F) \to (V, F)$ and $f:(V, F_1) \to (V, F_2)$ are both digraph isomorphisms.
\end{prop}

\begin{proof}
Let $L_1, L_2:V \to \{1, 2, \dots, n\}$ be the labeling functions associated with $F \cup F_1$, $F \cup F_2$ respectively.
First assume $F \cup F_1$ and $F \cup F_2$ both induce $\pi$.  Let $f = {L_2}^{-1} \circ L_1$; this is a bijection $V \to V$.  
By equation (\ref{eq:defpi}), 
\begin{align*}
(v, w) \in F & \iff (L_1(v), L_1(w) ) \in E(\pi)\\
&  \iff ({L_2}^{-1} \circ L_1(x), {L_2}^{-1} \circ L_1(x) ) = (f(v), f(w)) \in F.
\end{align*}
Similarly, by equation  (\ref{eq:defpic}),
\begin{align*}
(x, y) \in F_1 & \iff (L_1(x), L_1(y) ) \in E(\pi)^c\\
&  \iff ({L_2}^{-1} \circ L_1(x), {L_2}^{-1} \circ L_1(y) ) = (f(x), f(y)) \in F_2.
\end{align*}
Conversely, assume $g:V \to V$ is a digraph isomorphism $(V, F) \to (V, F)$ and  $(V, F_1) \to (V, F_2)$.  Then $L_2 \circ g \circ L_1^{-1}:\{1, \dots, n \} \to \{1, \dots, n \}$ is a bijection.  Moreover, for all distinct $i, j \in \{1, \dots, n \}$,  equation (\ref{eq:defL}) gives 
\begin{align*}
i < j & \iff ( {L_1}^{-1}(i) ,  {L_1}^{-1} (j) ) \in F \cup F_1 \\
& \iff  ( g \circ {L_1}^{-1}(i) , g \circ {L_1}^{-1} (j) ) \in F \cup F_2 \\
& \iff  [L_2 \circ  g \circ {L_1}^{-1}](i)  <[ L_2 \circ g \circ {L_1}^{-1}] (j) .
\end{align*}
The only order-preserving bijection on a finite totally ordered set is the identity, implying  $L_1 = L_2 \circ g$. 
Combining this with the fact that $g$ is an isomorphism on $(V, F)$, we get 
\begin{align*}
\Big\{\big(L_1(v), L_1 (w)\big)  \mid (v,w) \in F\Big \} & = \Big\{\big(L_2 \circ g(v), L_2 \circ g(w)\big)  \mid (v,w) \in F\Big \} \\
& = \Big\{\big(L_2\circ g(v), L_2 \circ g(w)\big)  \mid (g(v),g(w)) \in F\Big \} \\
& = \Big\{\big(L_2(x), L_2 (y)\big)  \mid (x,y) \in F\Big \}.
\end{align*} 
Hence the permutation induced by $F \cup F_1$ has the same inversion set as (and thus is equal to) the permutation induced by $F \cup F_2$.
\end{proof}

Applying this to the case where $D$ is the null graph, recall that up to isomorphism, there is only one transitive orientation on $K_n$. Moreover,  any bijection on the vertices is also an isomorphism on $D$.
In this extreme case,  the vertices are indistinguishable in both $D$ and $G^c$.  The following theorem generalizes from indistinguishable vertices to indistinguishable submodules.

\begin{thm} \label{thm:countTO}
Let $D=(V,F)$ be a permutation digraph, with underlying undirected graph $G$.
Let $M$ be a degenerate $0$-node of $G$, with $\mathcal{C}= \{X_1, X_2, \dots, X_m\}$ being a set of children of $M$ that induce isomorphic directed subgraphs of $D$; let
$g_{ij}$ denote an isomorphism  $D | X_i \to D | X_j$.
For any $\rho \in S_m$,
define $f_\rho:V \to V$ by:
\[
f_\rho(v) = 
\begin{cases}
g_{i\rho(i)}(v), \quad &  \text{ if } v \in X_i \text{ for some } 1 \leq i \leq m,\\
v, & \text{ if } v \notin \mathcal{C}.
\end{cases}
\]
For any transitive orientation $F_1$ on $G^c$, define another orientation  ${F_2}$  by  
\[
\big(f_\rho(v), f_\rho(w)\big) \in F_2  \iff (v, w ) \in  F_1.
\]
Then  $F_2$  is also a transitive orientation on $G^c$. 
Moreover,  $F \cup F_1$ and $F \cup F_2$ induce the same permutation.
\end{thm}

\begin{proof}
First we show that  $f_\rho$ is a digraph isomorphism $D \to D$.  Since $f_\rho$ is the identity outside $\bigcup \mathcal{C}$, we need only consider arcs with at least one endvertex  in $\bigcup \mathcal{C}$.  Since $M$ is a degenerate 0-node of $G$,  no vertices in different children of $M$ are adjacent, so we have only the following two cases.
\begin{enumerate}
\item If $v, w \in X_i$, then use the fact that $f_\rho | _{X_i} = g_{i \rho(i)}$ is an isomorphism.
\item If $v \in X_i$ and $w \notin \bigcup \mathcal{C}$, then $w$ must belong to another module $M'$. If $v, w$ are adjacent in $G$,  then $M$ and $M'$ must either be adjacent children of some other node, or submodules of adjacent children of some other node.  In either case, by Lemma \ref{lem:childorient}, all edges between vertices in $M'$ and $M$ are oriented the same way in $F$. Thus  $(v, w) \in F$ if and only if  $\big(f_\rho(v), f_\rho(w)\big) = \big(f_\rho(v),w ) \in F$.
\end{enumerate}
Ignoring orientation,  $f_\rho$ is a graph isomorphism $G \to G$, and so also $G^c \to G^c$.  By construction, $f_\rho$ will be a digraph isomorphism $(V, F_1) \to (V, F_2)$. Since $F_1$ is transitive, $F_2$ must also be transitive.  Now apply 
Proposition~\ref{prop:Phinotinjective}.
\end{proof}

\begin{figure}[h] 
   \centering
   \includegraphics[width=3.7in]{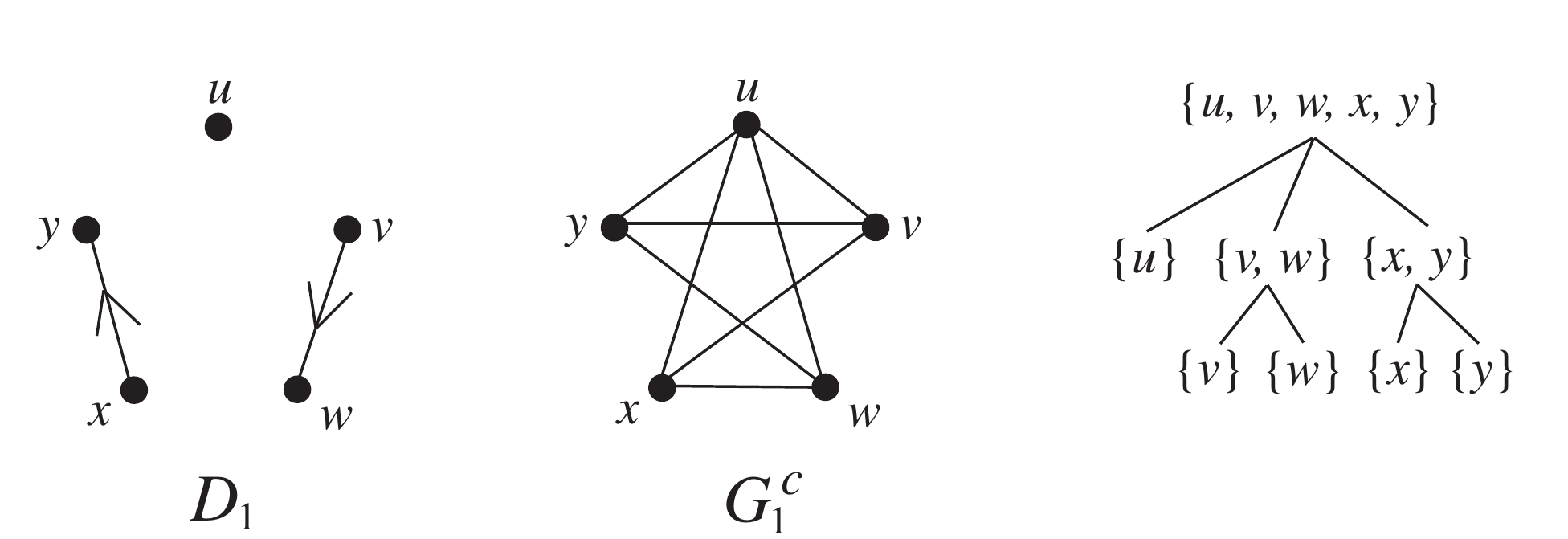} 
   \caption{Example~\ref{ex:MDTeg1} }
   \label{fig:countTO1}
\end{figure}

\begin{ex}\label{ex:MDTeg1}
 Figure~\ref{fig:countTO1} shows a permutation digraph $D_1$, the complement of the underlying undirected graph, $G_1^c$, and its modular decomposition tree.  
The root node $V$ is a degenerate 0-node of $G_1$ with $3$ children;  $\{v,w\}$ and $\{x,y\}$ are degenerate 1-nodes of $G_1$.  
Hence, $V$ is a degenerate 1-node of $G_1^c$ with $3$ children;  $\{v,w\}$ and $\{x,y\}$ are degenerate 0-nodes of $G_1^c$.  
By Corollary~\ref{cor:countTO}, the total number of different transitive orientations on $G_1^c$ is $3! = 6$.  However,  $\{v,w\}$ and $\{x,y\}$ induce isomorphic directed subgraphs of $D$,  so  by Theorem~\ref{thm:countTO}, orientations on $G_1^c$ that differ only by a permutation of these $2$ modules induce the same permutation.   Hence $D_1$ represents no more than $6/2!= 3$ permutations.  
\end{ex}

\begin{figure}[htbp] 
   \centering
   \includegraphics[width=3.7in]{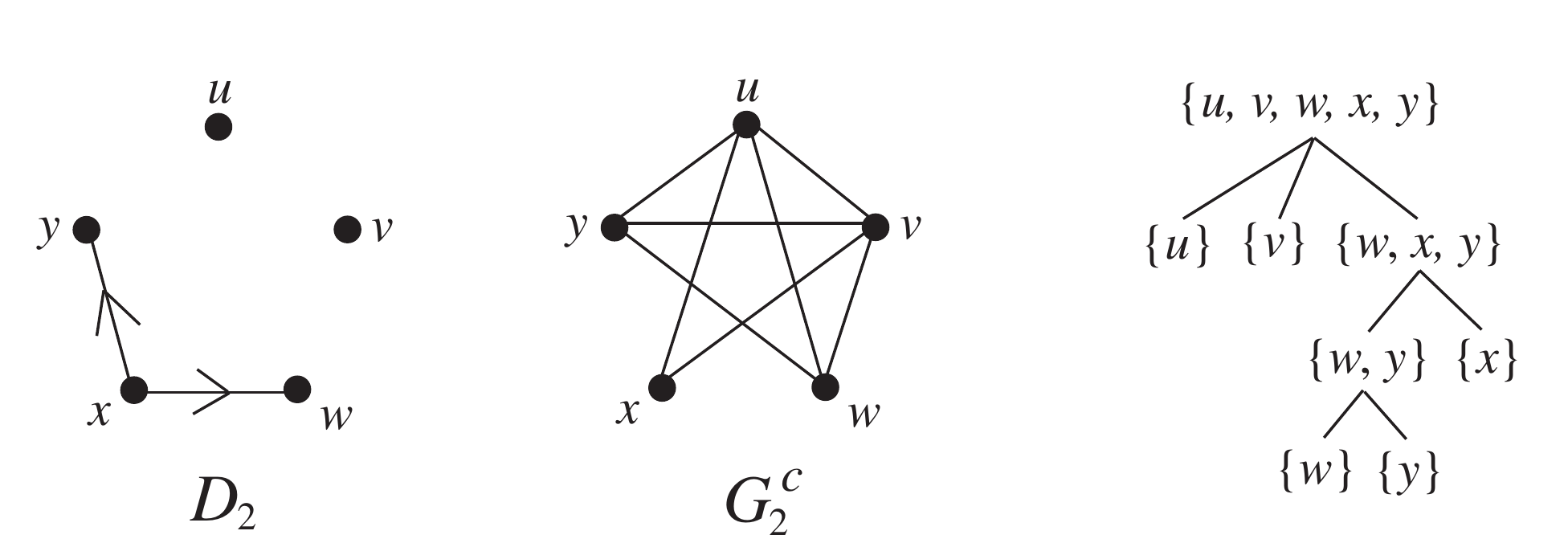} 
   \caption{Example~\ref{ex:MDTeg2} }
   \label{fig:countTO2}
\end{figure}

\begin{ex}\label{ex:MDTeg2}
In Figure~\ref{fig:countTO2}, $V$ is still a degenerate 1-node of $G_2^c$ with $3$ children;  $\{w, x, y\}$ is a degenerate 0-node and $\{w,y\}$ is a degenerate 1-node with 2 children.  By Corollary~\ref{cor:countTO}, the total number of different transitive orientations on $G_2^c$ is $3! \cdot 2!= 12$.  Since $\{u\}, \{v\}$ are isomorphic children of $V$, and $\{w\}, \{y\}$ are isomorphic children of $\{w,y\}$,  by Theorem~\ref{thm:countTO} (applied to both $D$ and $D |\{w,y\}$), $D_2$ represents $12/(2!\cdot 2!)= 3$ permutations.  
\end{ex}

Note that $D_1$ is isomorphic to its own reverse $-D_1$.  By Theorem~\ref{thm:permdigraphs}, the permutations represented by $D_1$ constitute one geo-equivalence class of $S_5$; in fact, it is class 2.4 in Table~\ref{table:S5}.  On the other hand, $D_2$ is not isomorphic to $-D_2$; the same argument as above shows that $-D_2$ also represents 3 permutations.  Hence the geo-equivalence class of permutations represented by either $D_2$ or $-D_2$ contains 6 permutations; it is class 2.3 in Table~\ref{table:S5}.
\medskip

A {\it cograph} is any graph whose modular decomposition tree contains no prime nodes.  As pointed out by Gallai,  such graphs and their complements are always transitively orientable, and hence all cographs are permutation graphs.  (The underlying undirected graphs in Examples \ref{ex:MDTeg1} and \ref{ex:MDTeg2} are both cographs.) As noted earlier, in \cite{Guru} Guruswami gives a generating function for the number of non-isomorphic cographs on $n$ vertices; in the same paper, he also shows that the number   of  $\pi \in S_n$ such that $G(\pi)$ is a cograph is $r_{n-1}$, where $(r_n)$ is the sequence of (large) Schr\"{o}der numbers (A006318 in the Online Encyclopedia of Integer Sequences \cite{oeis}).  Theorem~\ref{thm:countTO} allows us to determine the size of the geo-equivalence class of any such permutation.
\medskip

\begin{cor}\label{cor:cograph}
Let $\pi \in S_n$ such that $G(\pi) =G$ is a cograph, and let $D(\pi)= D$ be the corresponding permutation digraph.
Suppose the internal nodes of the modular decomposition tree of $G$ consist of:
\begin{itemize}
\item degenerate 1-nodes $N_1, \dots, N_t$;
\item degenerate 0-nodes $M_1, \dots, M_r$.
\end{itemize}
Suppose further that $M_i$ has $k_i$ children, which we divide up into isomorphism classes according to the directed subgraphs they induce in $D$: $n_{i1}$ children of isomorphism type $1$, $n_{i2}$ children of isomorphism type $2$, \dots and $n_{it}$ children of isomorphism type $t$.  Then the number of permutations represented by $D$ is 
\[
n(D) = \prod_{i = 1}^r \frac{k_i!} {n_{i1}! n_{i2}! \cdots n_{it}!}.
\]
The size of the geo-equivalence class $[\pi]$ is $2 n(D)$, unless $D \cong -D$, in which case it is $n(D)$.
\end{cor}

However,  the probability that a permutation has a cograph as its permutation graph approaches zero (the entry for  A00103 in OEIS \cite{oeis} gives the asymptotic behavior of the Schr\"{o}der numbers). Thus we now turn
 our attention to prime nodes. We begin with a lemma that gives a new perspective on the result of Corollary~\ref{cor:invrev}.

\begin{lem}\label{lem:invrev}
Let $G$ be a permutation graph and let $F, F_1$ be transitive orientations on $G, G^c$ respectively.  If $F \cup F_1$ induces  $\pi$, then:
\begin{enumerate}
\item \label{it1}$-F \cup F_1$ induces $\pi^{-1}$;
\item \label{it2}$F \cup -F_1$ induces $((\pi^c)^{-1})^c$;
\item \label{it3} $-F \cup -F_1$ induces $(((\pi^c)^{-1})^c)^{-1}$.
\end{enumerate}
\end{lem}
 
 \begin{proof}
 Let $L$ be the unique vertex labeling function corresponding to $F \cup F_1$.
For part~\ref{it1}, we first show that $\pi^{-1} \circ L$ is the unique vertex labeling function corresponding to $(-F) \cup F_1$; that is,  
\[
\pi^{-1} \circ L(v) < \pi^{-1} \circ L(w) \iff (v, w) \in -F \cup F_1.
\] 
Assume $\pi^{-1} \circ L(v) < \pi^{-1} \circ L(w)$.  
\begin{itemize}
\item If $L(v) > L(w)$, then by definition $(L(w), L(v))\in E(\pi)$, and so by (\ref{eq:defpi}),
$(w,v) \in F$ and hence $(v,w) \in -F$.
\item If $L(v) < L(w)$, then $(L(v), L(w)) \notin E(\pi)$, so  $(L(v), L(w)) \in E(\pi^c)$ and thus by (\ref{eq:defpic}), $(v,w) \in F_1$.
\end{itemize}
Conversely, assume $(v,w) \in -F \cup F_1$.
\begin{itemize}
\item If $(v, w) \in -F$, then $(w,v) \in F$ and so by (\ref{eq:defpi}), $(L(w), L(v)) \in E(\pi)$.  By Corollary~\ref{cor:inverse}, 
$
\big( \pi^{-1} \circ L (v), \, \pi^{-1} \circ L (w) \big) \in E(\pi^{-1}).
$
For this to be an inversion, it must be that  $\pi^{-1} \circ L(v) < \pi^{-1} \circ L(w)$.
\item If $(v,w) \in F_1$, then by (\ref{eq:defpic}), $(L(v), L(w)) \in E(\pi^c) = E(\pi)^c$.  Since  $(L(v), L(w))$ is not an inversion of $\pi$,  $\pi^{-1} \circ L(v) < \pi^{-1} \circ L(w)$.
\end{itemize}
The permutation defined by $\pi^{-1} \circ L$ has inversion set
\begin{align*}
\Big\{\big( &\pi^{-1} \circ L(v), \pi^{-1} \circ L(w)\big)   \mid (v,w) \in - F\Big \} \\
& = \Big\{\big(\pi^{-1} \circ L(v), \pi^{-1} \circ L(w)\big)  \mid (w,v) \in F\Big \} \\
& = \Big\{\big(\pi^{-1} ( L(v)), \pi^{-1} (L(w)) \big)  \mid (L(w),L(v)) \in E(\pi) \Big \}\\
& = \pi^{-1} \ast E(\pi) = E(\pi^{-1}),
\end{align*}
since the action of $\pi^{-1}$ on $E(\pi)$ is order-reversing.
\medskip

We can prove part 2 similarly by showing that  $(\pi^c)^{-1} \circ L$ is the unique vertex labeling function corresponding to $F \cup (-F_1)$.  To show that the corresponding inversion set is that of $((\pi^c)^{-1})^c$, recall (from Corollary~\ref{cor:reverse}) that the action of $(\pi^c)^{-1}$ on $E(\pi)$ is order-preserving.  Part 3 follows from parts 1 and 2.
\end{proof}

\begin{prop}\label{prop:primegraphs}
Let $\pi \in S_n$ such that $G(\pi) = G = (V, E)$ has only trivial modules, with $V$ being a prime node in the modular decomposition of $G$.  Then 
$[\pi]$ 
is a multiset of the form
\[
 \Big\{ \pi, \pi^{-1},  ((\pi^c)^{-1} )^c, ( ( (\pi^c)^{-1} )^c )^{-1} \Big\},
\]
which may contain four, two or one distinct permutation(s).
Moreover, this is the only geo-equivalence class of $S_n$  represented by a transitive orientation of $G$.
\end{prop}

\begin{figure}[h] 
   \centering
   \includegraphics[width=2.8in]{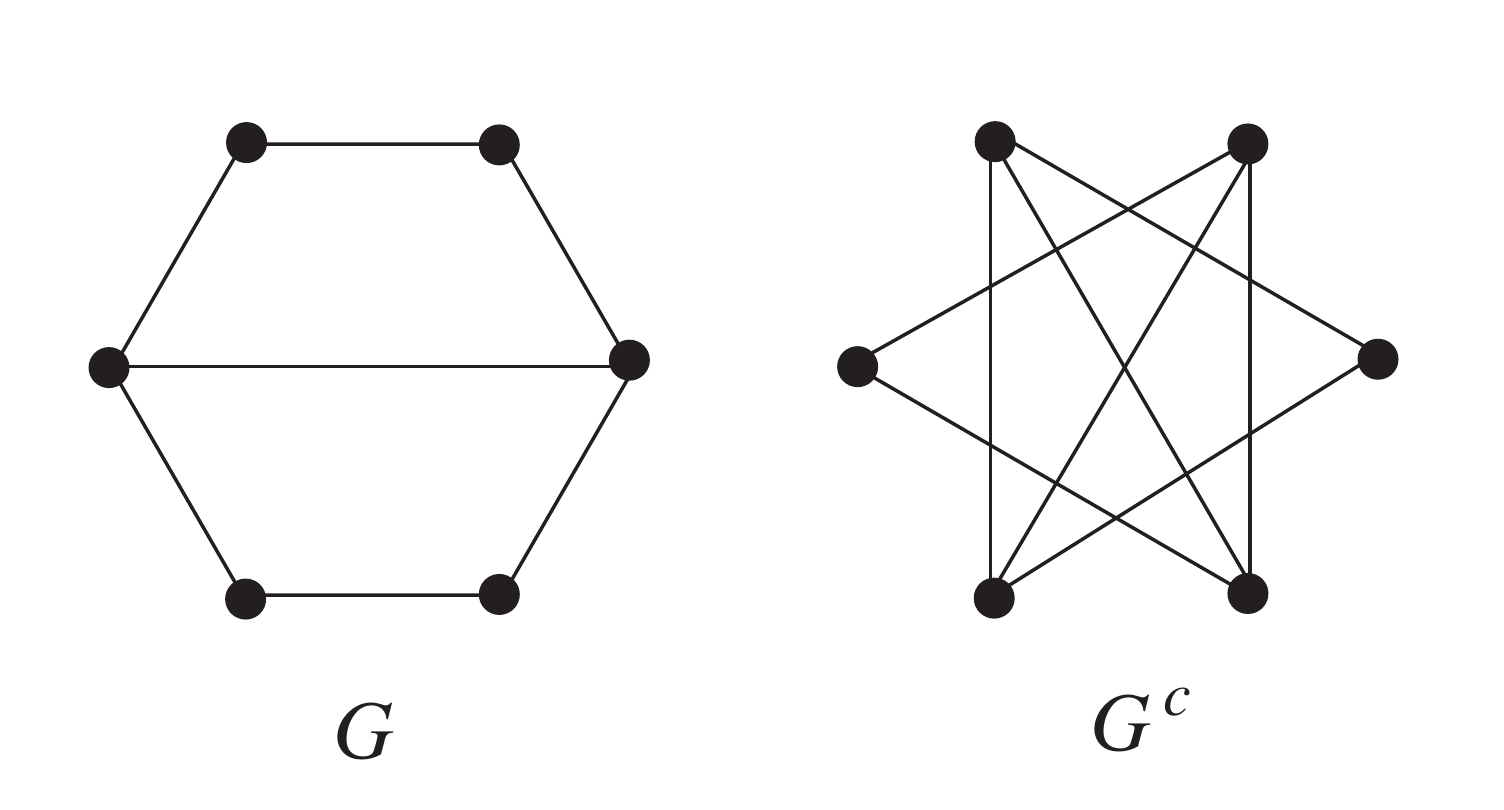} 
   \caption{A permutation graph representing  $[351624] = \{351624\}$.}
   \label{fig:allisone}
\end{figure}

\begin{proof}
Let $D(\pi) = (V, F)$ and $D(\pi^c) = (V, F_1)$.
By Proposition~\ref{prop:prime}, both $G$ and $G^c$ are uniquely orientable, so $F, -F$ and $F_1, -F_1$ are the only two transitive orientations on $G$ and $G^c$ respectively.
By Lemma~\ref{lem:invrev},  $D$ represents only $\pi$ and $((\pi^c)^{-1})^c$ and $-D$ represents only $\pi^{-1}$ and  $(( (\pi^c)^{-1} )^c )^{-1}$. 
\medskip

It is possible that $D \cong -D$; in that case, there exists a bijection $f:V \to V$ such that
$
(u,v) \in F \iff \big(f(u), f(v) \big) \in -F \iff \big(f(v), f(u)\big) \in F.
$
Now $f$ is also a graph isomorphism on $G^c$ and so $\{(f(x), f(y) ) | (x,y) \in F_1\}$ constitutes a transitive orientation on $G^c$.  Since $G^c$ is uniquely orientable,  it can only be either $F_1$ or $-F_1$.  
In the first case, $f:(V, F \cup F_1) \cong (V, -F \cup F_1)$; reversing orientations everywhere, we also get $f: (V, -F \cup -F_1) \cong (V, F \cup -F_1)$.  Lemma~\ref{lem:invrev} implies that 
$
\pi = \pi^{-1} \text {  and  }  ((\pi^c) ^{-1})^c)^{-1} = (\pi^c) ^{-1})^c.
$  
In the second case, $f:(V, F \cup F_1) \cong (V, -F \cup -F_1)$ and so by Lemma~\ref{lem:invrev}, $\pi = (((\pi^c) ^{-1})^c)^{-1}$, which (purely algebraically) implies 
$\pi^{-1} = ((\pi^c)^{-1})^c$.	
 Figure~\ref{fig:uniqor} shows one example of each case; for both examples, 
 $f(u)=u, f(v)=y, f(w)=x, f(x)=w$ and  $f(y) = v$.
 \medskip

\begin{figure}[htbp] 
   \centering
   \includegraphics[width=4.5in]{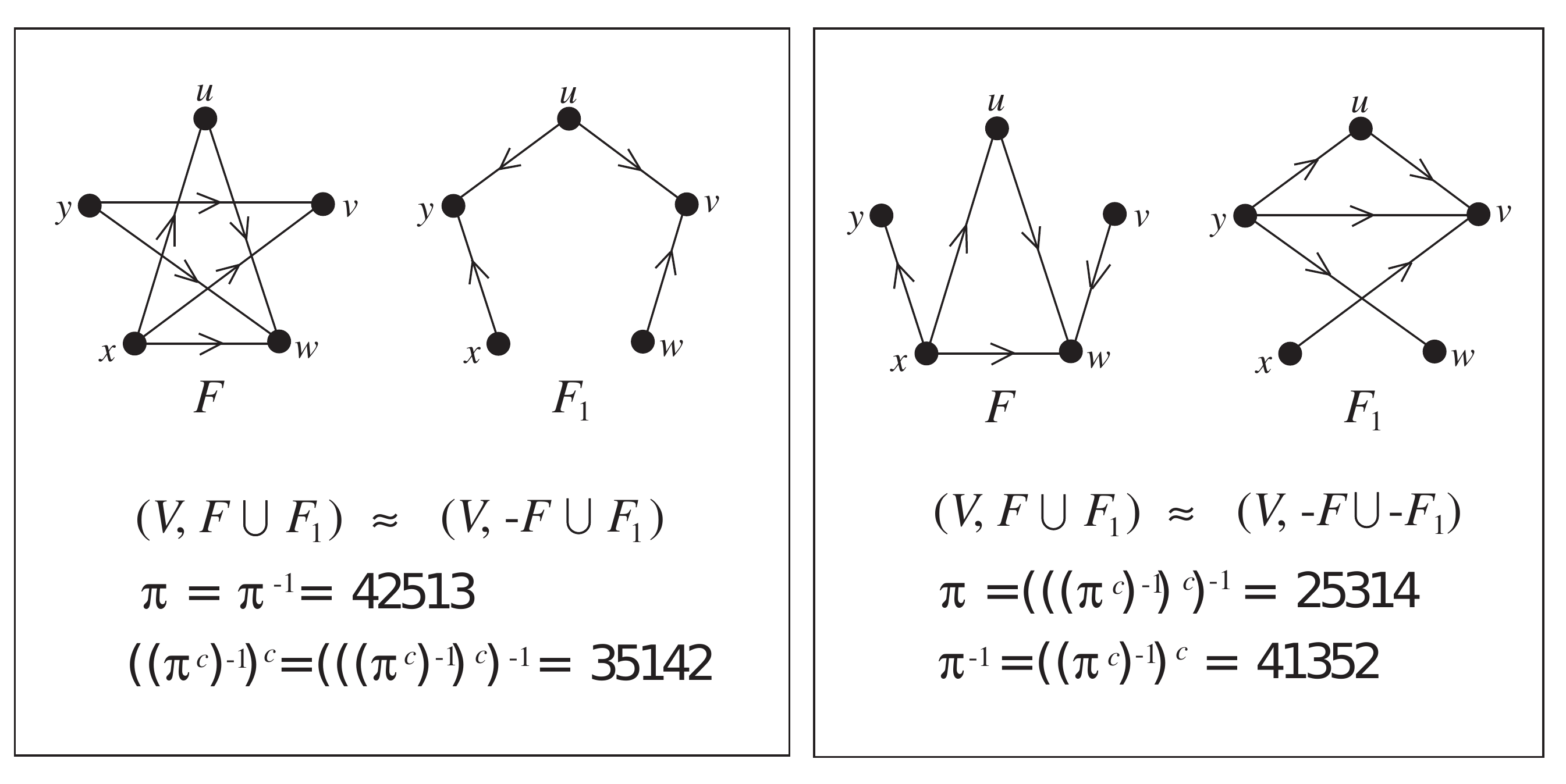} 
   \caption{Permutation graphs representing exactly two permutations.}
   \label{fig:uniqor}
\end{figure}

If $D \not \cong -D$ but $(V,F_1) \cong (V, -F_1)$, then an analogous argument shows that either 
$
\pi =( (\pi^c)^{-1})^c$ and  $ \pi^{-1} = (((\pi^c) ^{-1})^c)^{-1},
$
or 
$
\pi = (((\pi^c)^{-1})^c)^{-1}$  and  $ \pi^{-1} = ((\pi^c) ^{-1})^c.
$
\medskip

A third possibility is that 
\[
(V, F\cup F_1) \cong (V, - F \cup F_1) \cong (V, F\cup - F_1) \cong (V, - F\cup -F_1),
\]
in which case all four permutations are equal.  We leave it to the reader to verify that this situation occurs with the permutation graph in Figure~\ref{fig:allisone}.

\end{proof}  

\medskip

We now consider the situation where a prime node is one of several internal nodes in the modular decomposition tree.
By Corollary~\ref{cor:countTO}, each prime node contributes a factor of $2$ to the number of different transitive orientations on $G^c$.  However, it may contribute only a factor of 1 to the number of unrelated transitive orientations, as in the rather detailed special case described  below.

\begin{thm} \label{thm:primecountTO}
Let $D=(V,F)$ be a permutation digraph, with underlying undirected graph $G$.
Let $P$ be a prime node of $G$, with children $V^P= \{X_1, X_2, \dots, X_p\}$.
Let $Q(P), Q^c(P)$ denote the corresponding quotient graphs of $G, G^c$ respectively. 
Let $F^P, -F^P$ denote the (only) two transitive orientations on $Q(P)$ and let ${F_1}^P, -{F_1}^P$ denote the two transitive orientations on $Q^c(P)$.
Assume that there exists a bijection $f^P:V^P \to V^P$ such that
:\begin{itemize}
\item 
$f^P:(V^P, F^P) \to (V^P, F^P)$ and $f^P:(V^P, {F_1}^P) \to (V^P, -{F_1}^P)$ are both digraph isomorphisms, and 
\item 
 for all $1 \leq i \leq p$,   $X_i$ and $f^P(X_i)$ induce isomorphic subgraphs of $D$.
 \end{itemize}
 Let $g_i: D \mid X_i  \to D \mid f^P(X_i)$ be a digraph isomorphism.
Define $f:V \to V$ by:
\[
f(v) = 
\begin{cases}
g_i(v), \quad &  \text{ if } v \in X_i \text{ for some } 1 \leq i \leq p,\\
v, & \text{ otherwise.}
\end{cases}
\]
For any transitive orientation $F_1$ on $G^c$, define another orientation  ${F_2}$  by  
\[
\big(f(v), f(w)\big) \in F_2  \iff (v, w ) \in  F_1.
\]
Then  $F_2$  is also a transitive orientation on $G^c$. 
Moreover, 
  $F \cup F_1$ and $F \cup F_2$ induce the same permutation.
\end{thm}
  \begin{proof}
First we show that  $f$ is a digraph isomorphism $D \to D$.  Again, we only consider arcs with at least one endvertex  in $P$, but there are now three cases.
\begin{enumerate}
\item If $v, w \in X_i$, then use the fact that $f | _{X_i} = g_i$ is an isomorphism.
\item If $v \in X_i, \, w \in X_j$, then  use the assumption that $f^P$ is an isomorphism on $(V^P, F^P)$ and Lemma~\ref{lem:childorient}.
\item If $v \in X_i$ and $w \notin P$, then 
 use the same argument as in the proof of Theorem~\ref{thm:countTO}.
\end{enumerate}
The rest of the proof is exactly the same as that of Theorem~\ref{thm:countTO}.
\medskip

\begin{ex}
In Figure~\ref{fig:primeintEG},  the modular decomposition tree of $G^c$ has the root node $V=P$ as a prime node, $X_1=\{s, t\}$ as a degenerate  $0$-node and $X_2=\{w,x\}$, $X_3 =\{y,z\}$ as degenerate $1$-nodes. Thus by Corollary~\ref{cor:countTO}, there are $2 \cdot 2! \cdot 2! = 8$ different transitive orientations on $G^c$, but we can show that they all induce the same permutation represented by $D$.  First,  orientations on $G^c$ that differ only by the orientations on $X_2$ and $X_3$  induce the same permutation by Theorem~\ref{thm:countTO}. Next, define a bijection on the set of children of the prime node $V^P$ by
\[
f^P(X_1)=X_1,\,  f^P(\{u\})=\{v\}, \, f^P(\{v\})=\{u\},\,  f^P(X_2) = X_3 ,
\, f^P(X_3)=X_2.
\]
It is easy to verify that $f^P$ is a digraph isomorphism both $(V^P, F^P) \to (V^P, F^P)$ and $(V^P, {F_1}^P) \to (V^P, -{F_1}^P)$.
Let $g_1$ be the identity on $X_1$, and let $g_2:X_2 \to X_3$ 
and $g_3:X_3 \to X_2$ be 
any bijections. Constructing $f:V \to V$ as in Theorem~\ref{thm:primecountTO} (with $f(u) = v$ and $f(v)=u$),
we conclude that any two orientations on $G^c$ that differ only in that the orientation on $Q^c(P)$ is reversed induce the same permutation. (Note that in this example, $D \not \cong -D$, so the geo-equivalence class $[\pi]$ contains exactly two permutations, namely $\pi$ and $\pi^{-1}$.)
\end{ex}

\begin{figure}[htbp] 
   \centering
   \includegraphics[width=5.5in]{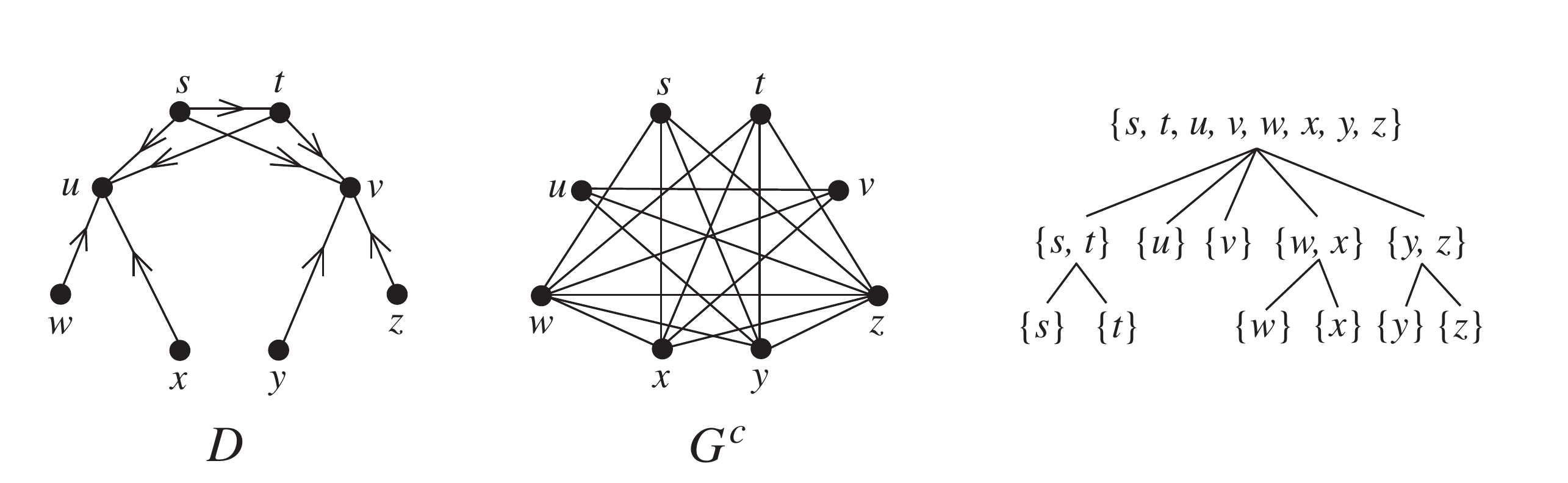} 
   \caption{A permutation digraph $D$ representing only $\pi = 51284367$.}
   \label{fig:primeintEG}
\end{figure}
\end{proof}




\section{Poset Structure}

Recall that for $\sigma, \pi \in S_n$,
\[
[\sigma] \preceq [\pi] \text{ in } [\mathcal{S}_n] \iff \Ksig \preceq \Kpi \text{ in }\mathcal{K}_{2,n},
\]
which holds if and only if there exists a geo-homomorphism $f:\Ksig \to \Kpi$ whose underlying map is a graph isomorphism. 
For strict precedence, $\Ksig$ must have strictly fewer edge crossings than $\Kpi$; equivalently, by Theorem~\ref{thm:inversion}, $|E(\sigma)| < |E(\pi)|$.  The proof of Theorem \ref{thm:complete} can easily be modified to yield the following.

 \begin{prop}\label{prop:homsm}
Let $\sigma, \pi \in S_n$.  Then $[\sigma] \prec [\pi]$ if and only if 
there exists $\rho \in S_n$ such that  
\begin{enumerate}
\item $\rho \ast E(\sigma) \subset   E(\pi) $;
\item $\rho$ is either order-preserving on  $E(\sigma)$ or order-reversing on $E(\sigma)$.
\end{enumerate}
\end{prop}

Of course, this result can be rephrased in terms of permutation digraphs.

\begin{prop}\label{prop:dihomsm}
Let $\sigma, \pi \in S_n$.  Then $[\sigma] \prec [\pi]$ if and only if $D(\sigma)$ is isomorphic to a proper directed subgraph of either $D(\pi)$ or  $D(\pi^{-1})$.
\end{prop}

\begin{cor} For all $n$, $ [\mathcal{S}_n]$ is a bounded poset, with first element $[1 2 \cdots  n]$ and  last element  $[ n (n-1) (n-2) \cdots  1]$.  
\end{cor}

\begin{figure}[h] 
   \centering
   \includegraphics[width=2in]{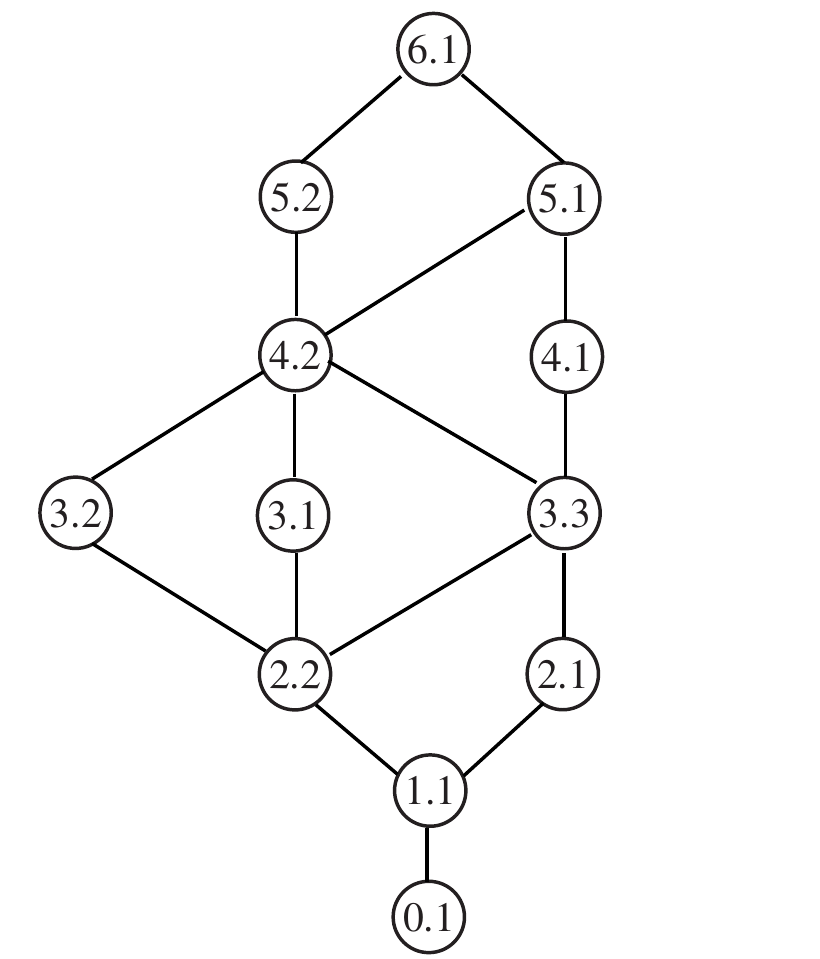} 
   \caption{The poset $\mathcal{S}_4$ ($\mathcal{K}_{2,4}$).}
   \label{fig:K24}
\end{figure}

\begin{figure}[h] 
   \centering
   \includegraphics[width=4.2in]{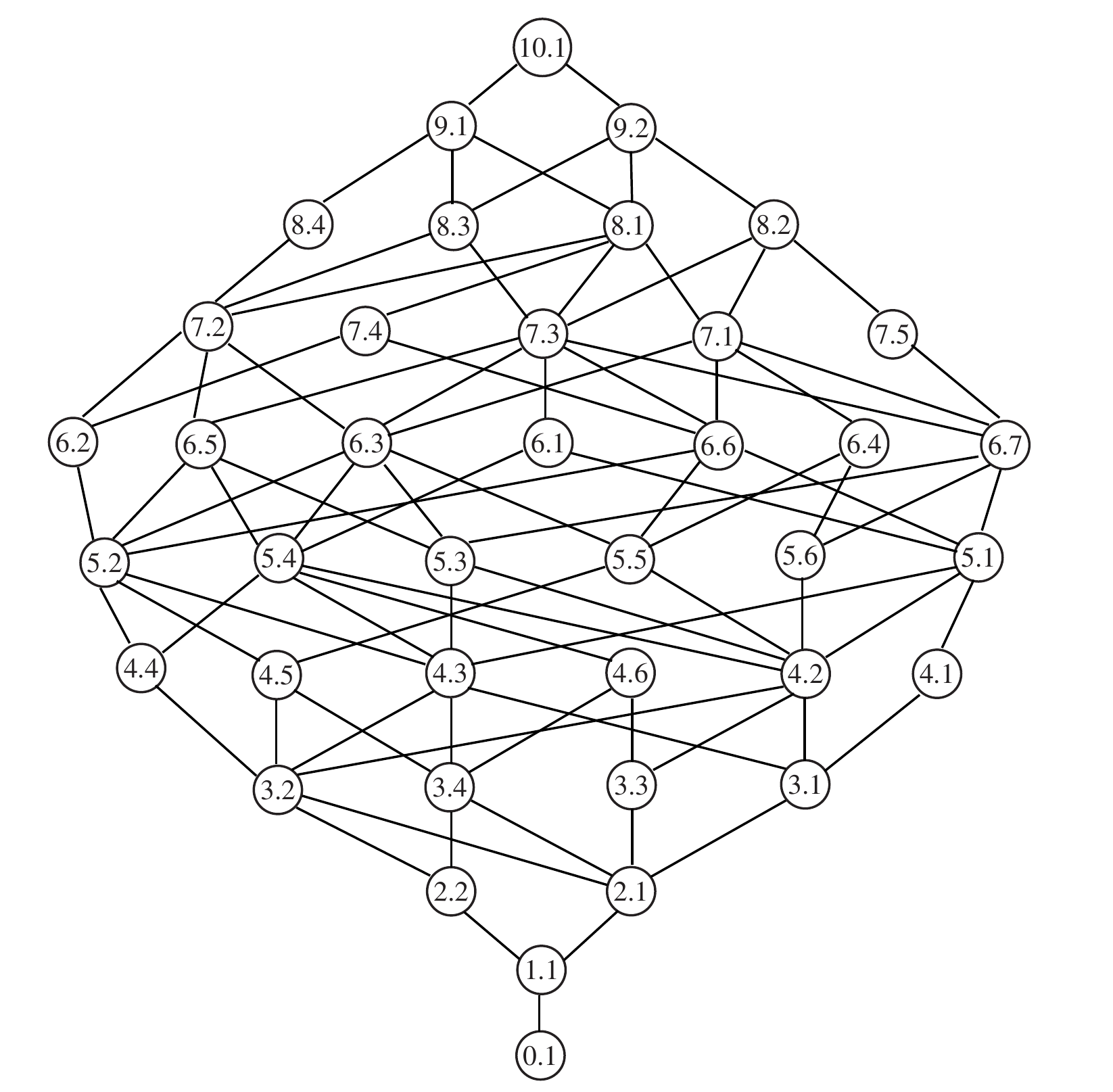} 
   \caption{The  poset $\mathcal{S}_5$ ($\mathcal{K}_{2,5}$).}
   \label{fig:S5poset}
\end{figure}

Proposition~\ref{prop:dihomsm} and visual inspection of the digraphs in Table~\ref{table:S4} determine the poset structure of $[\mathcal{S}_4]$; the Hasse diagram of this poset is given in Figure~\ref{fig:K24}.  Similarly, the industrious reader can fill in diagrams for the 39 geo-equivalence classes for $n=5$ in Table~\ref{table:S5} to obtain the poset structure of $[\mathcal{S}_5]$; the corresponding Hasse diagram is in Figure~\ref{fig:S5poset}.
\medskip

We compare this order on the geo-equivalence classes of $S_n$ to that induced by the weak left and right Bruhat orders, whose definitions we recall below. (For more on this order, see\cite{Bjorner}).  

\begin{defn}
Let $\sigma, \pi \in S_n$.  Then $\sigma$ strictly precedes $\pi$
\begin{itemize}
\item  in the {\it weak left Bruhat order} if and only if $E(\sigma) \subset E(\pi)$;
\item  in the {\it weak right Bruhat order} if and only if $E(\sigma^{-1}) \subset E(\pi^{-1})$.
\end{itemize}
\end{defn}

By Proposition~\ref{prop:dihomsm}, if $\sigma$ strictly precedes $\pi$ in either the weak left Bruhat order or the weak right Bruhat order, then 
$[\sigma] \prec [\pi]$ in $[\mathcal{S}_n]$.  In other words, the partial order on $[\mathcal{S}_n]$ is an extension of the order induced by the left and right weak Bruhat orders.

\medskip

 \begin{prop} \label{prop:adjtrans}\cite{Bjorner}
Let $\sigma, \pi \in S_n$ and let $\tau_i$ denote the adjacent transposition $ i \leftrightarrow i+1$ . Then $\pi$ covers $\sigma$ in
\begin{enumerate}
  \item the weak left Bruhat order $\iff$ $(\sigma(i), \sigma(i+1)) \not \in E(\sigma)$ and $\pi = \sigma \cdot \tau_i$ for some  $1\leq i <n$;
 \item the weak right Bruhat order $\iff$ $(i, i+1) \not \in E(\sigma)$ and $\pi = \tau_i \cdot \sigma$ for some  $1\leq i <n$.
 \end{enumerate}
 \end{prop}
 
 

\begin{ex}\label{ex:bruhat}
By Proposition~ \ref{prop:adjtrans}, $\sigma = 25314$ is covered in the weak left Bruhat order by
\[
 25314 \cdot \tau_1= 52314 \text{  and  }
 25314 \cdot \tau_4 = 25341,
 \]
 and in the weak right Bruhat order by
 \[
 \tau_2 \cdot 25314 = 35214 \text{  and  }
 \tau_3 \cdot 25314 = 25413.
 \]
From Table~\ref{table:S5}, $[\sigma]  = \{\sigma, \sigma^{-1}\} = \{ 25314, 41352\}$; by Proposition~ \ref{prop:adjtrans}, 
$41352$ is covered in the weak left Bruhat order by 
 \[
 41352 \cdot \tau_2 = 43152 \text{  and  }
 41352 \cdot \tau_3 = 41532.
\]
 and in the weak right Bruhat order by

 \[
\tau_1 \cdot 41352 = 42351 \text{  and  }
\tau_4 \cdot 41352 = 51342.
\]
 In terms of the class labels given in Table~\ref{table:S5}, the only covering relationships in $[\mathcal{S}_n]$ induced from the weak Bruhat orders are 
\[5.3 \prec 6.3\text{   and   } 5.3 \prec 6.7.
\]  
However, $\rho =23415$ is order-preserving on $E(\sigma)$
and
\[
\rho \ast E(\sigma) = \{(2,3), (2,4), (2,5), (4,5), (1,5)\}.
\]
This is not itself an inversion set (because its complement is not transitive), but it is a proper subset of $E( 35142)$.
Thus by Propostion~\ref{prop:homsm},  we also have $5.3 \prec 6.5$. This demonstrates that the partial order in $[\mathcal{S}_n]$  is a proper extension of that induced by the weak Bruhat orders.
  \end{ex}

\section{Open Questions}

\begin{enumerate}
\item Is there a (closed or recursive) formula for $a_n$, the number of geo-equivalence classes of $S_n$ (equivalently, the number of elements of $\mathcal{K}_{2,n}$)?

\item As shown in \cite{Bjorner},  $S_n$ is a graded lattice under the weak left Bruhat order, with the number of inversions serving as a rank function ({\it i.e.} if $\pi$ covers $\sigma$, then  $|E(\pi)| = |E(\sigma)| +1$).  Certainly $[\mathcal{S}_n]$ is not a lattice; Figure~\ref{fig:S5poset} shows that classes 8.3 and 8.1 both have classes 9.1 and 9.2 as suprema (and classes 9.1 and 9.2 have both 8.1 and 8.3 as infima).  However, $[\mathcal{S}_n] $ is a graded poset for $n \leq 5$, with the number of inversions as a rank function.  Rephrasing this using Theorem~\ref{thm:inversion}, the number of edge crossings serves as a rank function in the homomorphism poset $\mathcal{K}_{2,n}$, for $n \leq 5$.  In \cite{BCDM}, Boutin, Cockburn, Dean and Margea show that the homomorphism posets for paths $\mathcal{P}_n$, cycles $\mathcal{C}_n$ and cliques $\mathcal{K}_n$ are graded posets with the number of edge crossings as rank function for $n \leq 5$, but not for $n=6$.  In fact, for all $n \geq 6$, $\mathcal{P}_n$ and $\mathcal{C}_n$ are not graded posets.  Is $[\mathcal{S}_n]$ a graded poset for all $n$?
\end{enumerate}

\textbf{Acknowledgements.}  I am indebted to several people for their help with this paper. Debra Boutin and Alice Dean contributed to the proof that every realization of $K_{2,n}$ is geo-isomorphic to $\Kpi$ for some $\pi \in S_n$. At the University of Victoria, Peter Dukes  and Dennis Eppel both provided feedback and support. Thanks also to Rick Decker at Hamilton College for writing C++ code to  compute the number of geo-equivalence classes for $n = 7, 8$ and $9$.

\bibliographystyle{plain}
\bibliography{permutations(2)}{}

\end{document}